\documentclass[a4paper,12pt]{article}
    \usepackage[top=2cm,bottom=2cm,left=2.5cm,right=2.5cm]{geometry}
    \usepackage{cite, amsmath, amssymb}
    \usepackage[margin=0.5cm,%
                font=small,%
                format=hang,%
                labelsep=period,%
                labelfont=bf]{caption}
\usepackage{amsmath}
\usepackage{amsfonts}
\usepackage{amssymb}
\usepackage{cite}
\usepackage{amsthm}
\usepackage{makeidx}
\usepackage{graphicx}
\usepackage{mathrsfs,xcolor,multicol}
\usepackage{enumerate}
\usepackage{blkarray}
\usepackage{bm}
\usepackage{setspace}
\usepackage{float}
\usepackage{authblk}
\usepackage{multirow}
\usepackage{booktabs}
\usepackage[ruled,vlined]{algorithm2e}
\usepackage{blkarray}
\usepackage{authblk}
\usepackage{multirow}
\usepackage{lineno}
\usepackage{lscape}
\usepackage{blkarray, bigstrut}
\usepackage{authblk}
\usepackage{multirow}

\usepackage{fancyhdr}
\usepackage{lastpage}

\usepackage{amsmath, amssymb, booktabs, lastpage,}
\usepackage{multicol}
\usepackage{hyperref}
\usepackage{blkarray}
\usepackage{empheq} 
\usepackage{latexsym} 
\usepackage{array}
\usepackage{tikz-cd}
\usepackage{tikz}
\usetikzlibrary{arrows}
\usetikzlibrary{automata,positioning}
\usetikzlibrary{circuits.logic.US,circuits.logic.IEC,fit}
\newcommand*\circled[1]{\tikz[baseline=(char.base)]{
            \node[shape=circle,draw,inner sep=2pt] (char) {#1};}}

\numberwithin{table}{section}
\numberwithin{equation}{section}

\newcommand{\matindex}[1]{\mbox{\scriptsize#1}}

\theoremstyle{plain}
\newtheorem{theorem}{Theorem}[section]
\newtheorem{proposition}[theorem]{Proposition}
\newtheorem{definition}[theorem]{Definition}

\newtheorem{remark}[theorem]{Remark}

\newcommand{\NN}{\mathcal{N}}

\usepackage{authblk}
\author[1]{ \textbf{Al Jay Lan J. Alamin}}
\author[1]{ \textbf{Melquezedec James T. Cruz}}
\author[1]{ \textbf{Bryan S. Hernandez}}
\author[2,3,4]{ \textbf{Eduardo R. Mendoza}}

\affil[1]{\small \textit{Institute of Mathematics, University of the Philippines Diliman, Quezon City 1101, Philippines}}
\affil[2]{\small \textit{Systems and Computational Biology Research Unit, Center for Natural Sciences and Environmental Research, Manila 0922, Philippines}}
\affil[3]{\small \textit{Department of Mathematics and Statistics, De La Salle University, Taft Avenue, Manila, 0922, Philippines}}
\affil[4]{\small \textit{Max Planck Institute of Biochemistry, Martinsried near Munich, 82152, Germany}}
\affil[*]{Email addresses: \texttt{ajalamin@up.edu.ph},
\texttt{mtcruz16@up.edu.ph}, \texttt{bshernandez@up.edu.ph} (corresponding author), \texttt{eduardo.mendoza@dlsu.edu.ph}}

\title{\textbf{The Long-Term Impact of Direct Capture Approaches to Carbon Dioxide Removal}}
\date{}

\begin{document}
\maketitle
\begin{abstract} 
Understanding the similarities and differences of the long term impact of different carbon dioxide removal (CDR) techniques is essential in determining the most effective and sustainable strategies to mitigate climate change. In particular, direct ocean capture (DOC) has emerged as a promising approach. In contrast to direct air capture (DAC) which separates carbon dioxide from the atmosphere, DOC performs the separation directly from seawater before storing it in geological reservoirs. In this study, we construct and analyze a kinetic system for CDR via DOC using chemical reaction network theory. Our analysis reveals the necessary conditions for the existence of positive steady states and highlights the potential for multistationarity, where the carbon cycle may admit multiple positive steady states, emphasizing the critical importance of addressing tipping points, thresholds beyond which the system could undergo irreversible changes. Furthermore, we examine conditions under which certain carbon pools exhibit absolute concentration robustness, remaining resistant to change regardless of initial conditions. We also determine the conditions for the carbon reduction capability of the model with the DOC intervention.
Importantly, a comparative analysis is then presented, where we compare the DOC model with the well-established DAC model by Fortun et al., and explore an integrated DOC-DAC approach for CDR. This comparison is important given that DAC is already being implemented in large-scale projects, while DOC remains in its early stages with limited trials and is geographically constrained to oceanic vicinity. Our comparative modeling framework provides valuable insights into the long-term impacts and complementary roles of DOC, DAC, and their integration into broader CDR strategies for climate mitigation.
	
\end{abstract}

\thispagestyle{empty}

\section{Introduction}
\label{intro}

Over the past several hundred years, the expansion of society's consumption of fossil fuels and extensive alteration of the terrestrial biosphere has led to a dramatic rise in levels of carbon dioxide and other greenhouse gases in the atmosphere. The resulting climate change is one of the most serious issues society is facing today. It is challenging to significantly cut down on $\text{CO}_2$ emissions since this modern world relies heavily on fossil fuels to keep economy running \cite{nationalacademies2022research}.

The Earth's carbon cycle is a complex and dynamic system that plays an important role in regulating the climate of our planet and sustaining life. It involves the exchange of carbon between terrestrial ecosystems, the atmosphere, and the oceans. Understanding the intricacies of this cycle is important for predicting the impacts of activities, such as anthropogenic carbon dioxide (CO$_2$) emissions in the atmosphere, which affect global climate change and for developing strategies to mitigate these effects \cite{hansen2016assessing,rosenzweig2008attributing,falkowski2000global}.

So far, efforts to remove excess $\text{CO}_2$ from the air have largely focused on what can be done on the land, such as growing trees or building direct air capture plants \cite{lebling2022ocean, rudee2020restoring, lebling20226things}. However, a growing number of researchers, companies and even national governments have begun to look at the ocean as a potential location for carbon dioxide removal \cite{nationalacademies2022research, stripe2021stripe}.

In the fight against climate change, Carbon Dioxide Removal (CDR) technologies are essential for reducing CO$_2$ levels.
We know that the ocean is good at sequestering carbon because it has already absorbed 30\% of the $\text{CO}_2$ - and 90\% of excess heat - caused by human activities, significantly dampening the impacts of climate change \cite{gruber2019oceanic, portner2019ipcc}. In total, the ocean holds around 42 times more carbon than the atmosphere \cite{lebling2022ocean, friedlingstein2022essd}.

CDR through Direct Ocean Capture (DOC) incorporates novel electrochemical engineering techniques where dissolved CO$_2$ is separated from seawater and stored in geological stock, with a physical process similar to Direct Air Capture (DAC) \cite{fortun2024determining}.

A key to understanding the DOC system is the application of chemical reaction network theory (CRNT). CRNT is particularly valuable for analyzing the structural and dynamical behavior of a system with uncertain or variable parameters. In particular, we explore crucial properties of DOC systems using CRNT: existence of positive steady states, \emph{multistationarity} and \emph{absolute concentration robustness (ACR)}. We also identify conditions for the carbon reduction capability of the DOC system.

Studying the steady states of a system provides us an understanding of its long-term behavior and helps us determine its stability.
Furthermore, understanding the complexities of climate change requires a thorough examination of climate tipping points. These points denote critical thresholds where the climate system undergoes changes that could lead to irreversible impacts. Predicting and comprehending these tipping points is crucial for developing effective strategies to mitigate the impacts of climate change.  \emph{Multistationarity}, associated to tipping points, describes how a system could swiftly and irreversibly switch to another state. In the context of chemical reaction networks, multistationarity refers to the system's ability to maintain multiple steady states under identical parameters, including the same set of rate constants and conserved quantities. 
On the other hand, \emph{ACR} ensures the maintenance of the concentration level of key species despite changes in initial conditions. For DOC systems, achieving ACR is critical to maintaining robustness in carbon capture and storage processes over the long term.

We then conduct a comparative analysis of the structural and dynamic properties of the Direct Ocean Capture (DOC) model, alongside the well-established Direct Air Capture (DAC) model by Fortun et al. This comparison is crucial given that DAC is a well-established technology with large-scale projects already, while DOC is still in the early stages, with only a few trials conducted so far. Furthermore, DOC is geographically constrained to oceanic vicinity.

The integration of multiple technologies, i.e., the integrated DOC-DAC approach, is likely to be necessary for large-scale carbon reduction, and our study demonstrates how this can be effectively modeled within the CDR framework.

\section{Preliminaries}
\label{prelim}

\subsection{Chemical reaction networks}
\label{subs:CRN}

A \emph{chemical reaction network} or simply \emph{CRN} 
is a triple of nonempty finite sets, 
where
    \begin{itemize}
        \item[i.] $\mathcal{S} = \left\{A_1, A_2, \ldots, A_m \right\}$ is the set of \emph{species},
        \item[ii.] $\mathcal{C} = \{C_1,C_2, \ldots, C_n\}$ is the set of \emph{complexes} that are non-negative linear combinations of the species, and
        \item[iii.] $\mathcal{R} = \{R_1,R_2, \ldots, R_n\} \subset \mathcal{C} \times \mathcal{C}$ is the set of \emph{reactions}.
    \end{itemize}

A reaction $(C_i,C_j) \in \mathcal{R}$ is typically represented as $C_i \to C_j$. The complex $C_i$ is called the \emph{reactant complex} and $C_j$ is called the \emph{product complex}. The \emph{reaction vector} for this reaction is defined by the difference $C_j - C_i$.
Furthermore, the linear subspace $S$ of $\mathbb{R}^m$ spanned by the reaction vectors is called the \emph{stoichiometric subspace} of a given network, i.e., $S = \mathrm{span}\{C_j - C_i \in \mathbb{R}^m \mid C_i \rightarrow C_j \in \mathcal{R}\}$.

Consider the CRN, hereafter referred to as the DOC (direct ocean capture) network, which consists of the following seven reactions:
\allowdisplaybreaks
\begin{align*}
    &R_1: A_1 + 2A_2 \rightarrow 2A_1 + A_2\\
    &R_2: 2A_1 +A_2 \rightarrow A_1 + 2A_2\\
    &R_3: A_2 \rightarrow A_3\\
    &R_4: A_3 \rightarrow A_2\\
    &R_5: A_4 \rightarrow A_2\\
    &R_6: A_{17} \rightarrow A_4\\
    &R_7: A_3 \rightarrow A_{17}.
\end{align*}

The network has $m=5$ species ($A_1$, $A_2$, $A_3$, $A_4$, and $A_{17}$). Furthermore, it has $n=6$ complexes ($A_1+2A_2$, $2A_1+A_2$, $A_2$, $A_3$, $A_4$ and $A_{17}$) and has $r=7$ reactions ($R_1,R_2,\ldots,R_7$).

The {\emph{molecularity matrix}} $Y$, is an $m\times n$ matrix where $Y_{ij}$ is the stoichiometric coefficient of species $A_i$ in complex $C_j$. The {\emph{incidence matrix}} $I_a$ is an $n\times r$ matrix where 
$${\left( {{I_a}} \right)_{ij}} = \left\{ \begin{array}{rl}
 - 1&{\rm{ if \ }}{C_i}{\rm{ \ is \ in \ the\ reactant \ complex \ of \ reaction \ }}{R_j},\\
 1&{\rm{  if \ }}{C_i}{\rm{ \ is \ in \ the\ product \ complex \ of \ reaction \ }}{R_j},\\
0&{\rm{    otherwise}}.
\end{array} \right.$$
The {\emph{stoichiometric matrix}} $N$ is the $m\times r$ matrix given by 
$N=YI_a$.

The \emph{deficiency} of a CRN is $\delta=n-\ell-s$ where $n$ is the number of complexes, $\ell$ is the number of connected components, and $s$ is the rank of the stoichiometric matrix of the network.

For our network, the molecularity, incidence, and stoichiometric matrices are given by

\begin{center}
    $Y=\begin{blockarray}{cccccccc}
        & \matindex{$A_1+2A_2$} & \matindex{$2A_1+A_2$} & \matindex{$A_2$} & \matindex{$A_3$} & \matindex{$A_4$} & \matindex{$A_{17}$} \\
        \begin{block}{c[ccccccc]}
        \matindex{$A_1$} & 1 & 2 & 0 & 0 & 0 & 0 \\
        \matindex{$A_2$} & 2 & 1 & 1 & 0 & 0 & 0 \\
        \matindex{$A_3$} & 0 & 0 & 0 & 1 & 0 & 0 \\
        \matindex{$A_4$} & 0 & 0 & 0 & 0 & 1 & 0 \\
        \matindex{$A_{17}$} & 0 & 0 & 0 & 0 & 0 & 1 \\
        \end{block}
        \end{blockarray}$,
\end{center}
\begin{center}
    $I_a=\begin{blockarray}{cccccccc}
        & \matindex{$R_1$} & \matindex{$R_2$} & \matindex{$R_3$} & \matindex{$R_4$} & \matindex{$R_5$} & \matindex{$R_6$} & \matindex{$R_7$} \\
        \begin{block}{c[ccccccc]}
        \matindex{$A_1+2A_2$} & -1 & 1 & 0 & 0 & 0 & 0 & 0 \\
        \matindex{$2A_1+A_2$} & 1 & -1 & 0 & 0 & 0 & 0 & 0 \\
        \matindex{$A_2$} & 0 & 0 & -1 & 1 & 1 & 0 & 0 \\
        \matindex{$A_3$} & 0 & 0 & 1 & -1 & 0 & 0 & -1 \\
        \matindex{$A_4$} & 0 & 0 & 0 & 0 & -1 & 1 & 0 \\
        \matindex{$A_{17}$} & 0 & 0 & 0 & 0 & 0 & -1 & 1 \\
        \end{block}
        \end{blockarray}$,
\end{center}
and
\begin{center}
    $N=YI_a=\begin{blockarray}{cccccccc}
        & \matindex{$R_1$} & \matindex{$R_2$} & \matindex{$R_3$} & \matindex{$R_4$} & \matindex{$R_5$} & \matindex{$R_6$} & \matindex{$R_7$} \\
        \begin{block}{c[ccccccc]}
        \matindex{$A_1$} & 1 & -1 & 0 & 0 & 0 & 0 & 0 \\
        \matindex{$A_2$} & -1 & 1 & -1 & 1 & 1 & 0 & 0 \\
        \matindex{$A_3$} & 0 & 0 & 1 & -1 & 0 & 0 & -1 \\
        \matindex{$A_4$} & 0 & 0 & 0 & 0 & -1 & 1 & 0 \\
        \matindex{$A_{17}$} & 0 & 0 & 0 & 0 & 0 & -1 & 1 \\
        \end{block}
        \end{blockarray}$.
\end{center}
The deficiency of the DOC network is $\delta=n-\ell-s=6-2-4=0$ because there are six complexes, two connected components, and the rank of $N$ is four.

A CRN is weakly reversible if each of its reactions is contained in a directed cycle. Since each reaction in the DOC network belongs to a cycle, it is a weakly reversible network.

Therefore, the DOC network is a weakly reversible and deficiency zero network.

\subsection{Chemical kinetic systems}

A \emph{kinetics} for a reaction network $\mathcal{N}=(\mathcal{S}, \mathcal{C}, \mathcal{R})$ is an assignment to
each reaction $C_i \to C_j \in \mathcal{R}$ of a continuously differentiable {rate function} $\mathcal{K}_{C_i\to C_j}: \mathbb{R}^\mathcal{S}_{{\geq} 0} \to \mathbb{R}_{\ge 0}$ such that this positivity condition holds:
$\mathcal{K}_{C_i\to C_j}(c) > 0$ if and only if ${\sf{supp \ }} C_i \subset {\sf{supp \ }} c$, Here, ${\sf{supp \ }} C_i$ refers to the support of the vector $C_i$, which is the set of species with nonzero coefficient in $C_i$.
Hence, the pair $\left(\mathcal{N},\mathcal{K}\right)$ is called a \emph{chemical kinetic system}.

The \emph{species formation rate function} (SFRF) of $(\mathcal{N},\mathcal{K})$ is defined as $$f\left( x \right) = \displaystyle \sum\limits_{{C_i} \to {C_j} \in \mathcal{R}} {{\mathcal{K}_{{C_i} \to {C_j}}}\left( x \right)\left( {{C_j} - {C_i}} \right)}$$ with $x$ a vector of concentrations of the species that change over time. Equivalently,
$f(x) = N\mathcal{K}(x)$ where $N$ is the stoichiometric matrix of $\mathcal{N}$ and $\mathcal{K}(x)$ is the vector of rate functions.
The system of \emph{ordinary differential equations} (ODEs) of a chemical kinetic system is given by $\dfrac{{dx}}{{dt}} = f\left( x \right)$. A \emph{positive steady state} is a positive vector that makes each time derivative equal to zero. Thus, the set of positive steady states of a chemical kinetic system $\left(\mathcal{N},\mathcal{K}\right)$ is given by
${E_ + }\left(\mathcal{N},\mathcal{K}\right)= \left\{ {x \in \mathbb{R}^m_{>0}|f\left( x \right) = 0} \right\}.$

\subsection{Power law systems}

A \emph{power law kinetics} has the form ${\mathcal{K}_i}\left( x \right) = {k_i}\prod\limits_j {{x_j}^{{F_{ij}}}}$  
for each reaction $i =1,\ldots,r$ where ${k_i} \in {\mathbb{R}_{ > 0}}$ and ${F_{ij}} \in {\mathbb{R}}$. The $r \times m$ matrix  $F=\left[ F_{ij} \right]$ is called the \emph{kinetic order matrix} that contains the kinetic order values $F_{ij}$, and $k_i$ is called the $i$th rate constant. A \emph{power law system} is a CRN endowed with power law kinetics.

Specifically, if each kinetic order row contains the stoichiometric coefficients of each reactant for the associated reaction in the network, then the system follows the well-known \emph{mass action kinetics}.

\subsection{Network decomposition}
\label{subs:network:decomposition}
We can decompose a CRN into pieces of networks called \emph{subnetworks} by partitioning its reaction set into disjoint subsets.
A network decomposition $\NN = \NN_1 \cup \NN_2 \cup \ldots \cup \NN_k$ is said to be \emph{independent} if its stoichiometric subspace is a direct sum of the stoichiometric subspaces of its subnetworks. An equivalent condition is to show that the rank of the stoichiometric matrix of the whole network is the sum of the ranks of the stoichiometric matrices of its subnetworks.

This concept of independent decomposition is important to our study, as it establishes a significant relationship between the structure of the set of positive steady states of a given network and its independent subnetworks. The following result by M. Feinberg, which we call Feinberg Decomposition Theorem, highlights this relationship \cite[Appendix 6.A]{feinberg2019crnt}.
\begin{theorem}
    \label{thm:feinberg-decomp}
    Let $(\mathcal N, \mathcal K)$ be a chemical kinetic system. Suppose $\mathcal N$ is decomposed into $k$ subnetworks, say $\mathcal N_1, \mathcal N_2, \dots, \mathcal N_k$, and denote the restriction of $\mathcal K$ to the restrictions in $\mathcal N_i$ as $\mathcal K_i$. 
    If the network decomposition is independent, then \[\bigcap_{i=1}^k E_+(\mathcal N_i, \mathcal K_i) = E_+(\mathcal N, \mathcal K).\]
\end{theorem}

To get the finest independent decomposition (independent decomposition with maximum number of subnetworks), a MATLAB program was provided in \cite{LubeniaINDECS}. By entering the DOC network and applying the program, we obtain the following such decomposition: $\NN_1=\{R_1,R_2\}$ and $\NN_2=\{R_3,R_4,\ldots,R_7\}$.

Recall from Section \ref{subs:CRN} that the stoichiometric matrix of the whole DOC network is

    \begin{center}
    $N=\begin{blockarray}{cccccccc}
        \matindex{$R_1$} & \matindex{$R_2$} & \matindex{$R_3$} & \matindex{$R_4$} & \matindex{$R_5$} & \matindex{$R_6$} & \matindex{$R_7$}\\
        \begin{block}{[ccccccc]c}
        1 & -1 & 0 & 0 & 0 & 0 & 0 & \matindex{$A_1$}\\
        -1 & 1 & -1 & 1 & 1 & 0 & 0 & \matindex{$A_2$}\\
        0 & 0 & 1 & -1 & 0 & 0 & -1 & \matindex{$A_3$}\\
        0 & 0 & 0 & 0 & -1 & 1 & 0 & \matindex{$A_4$}\\
        0 & 0 & 0 & 0 & 0 & -1 & 1 & \matindex{$A_{17}$}\\
        \end{block}
        \end{blockarray}$.
    \end{center}

Furthermore, the stoichiometric matrices of the two subnetworks ($\NN_1$ and $\NN_2$) are
    $$N_1=\begin{blockarray}{ccc}
        \matindex{$R_1$} & \matindex{$R_2$}\\
        \begin{block}{[cc]c}
        1 & -1 & \matindex{$A_1$}\\
        -1 & 1 & \matindex{$A_2$}\\
        0 & 0 & \matindex{$A_3$}\\
        0 & 0 & \matindex{$A_4$}\\
        0 & 0 & \matindex{$A_{17}$}\\
        \end{block}
        \end{blockarray}
    {\text{and }}
    N_2=\begin{blockarray}{cccccc}
         \matindex{$R_3$} & \matindex{$R_4$} & \matindex{$R_5$} & \matindex{$R_6$} & \matindex{$R_7$}\\
        \begin{block}{[ccccc]c}
         0 & 0 & 0 & 0 & 0 & \matindex{$A_1$}\\
         -1 & 1 & 1 & 0 & 0 & \matindex{$A_2$}\\
         1 & -1 & 0 & 0 & -1 & \matindex{$A_3$}\\
         0 & 0 & -1 & 1 & 0 & \matindex{$A_4$}\\
         0 & 0 & 0 & -1 & 1 & \matindex{$A_{17}$}\\
        \end{block}
        \end{blockarray}.$$
    Since rank $N=4$, rank $N_1 = 1$ and rank $N_2 = 3$. Then, the sum of the ranks of the stoichiometric matrices of the subnetworks is the rank of the stoichiometric matrix of the whole network. Indeed, the decomposition is independent.

\section{Results and Discussion}
\label{results}

\subsection{The direct ocean capture system}

The Direct Ocean Capture (DOC) system is based on a three-compartment biochemical framework of the pre-industrial model of Anderies et al. \cite{anderies2013topology}, which describes the carbon cycle interactions through the transfer of carbon between the land ($A_1$), atmosphere ($A_2$), and ocean ($A_3$).

\begin{figure}[H]
    \begin{center}
    \includegraphics[width=11cm,height=6cm,keepaspectratio]{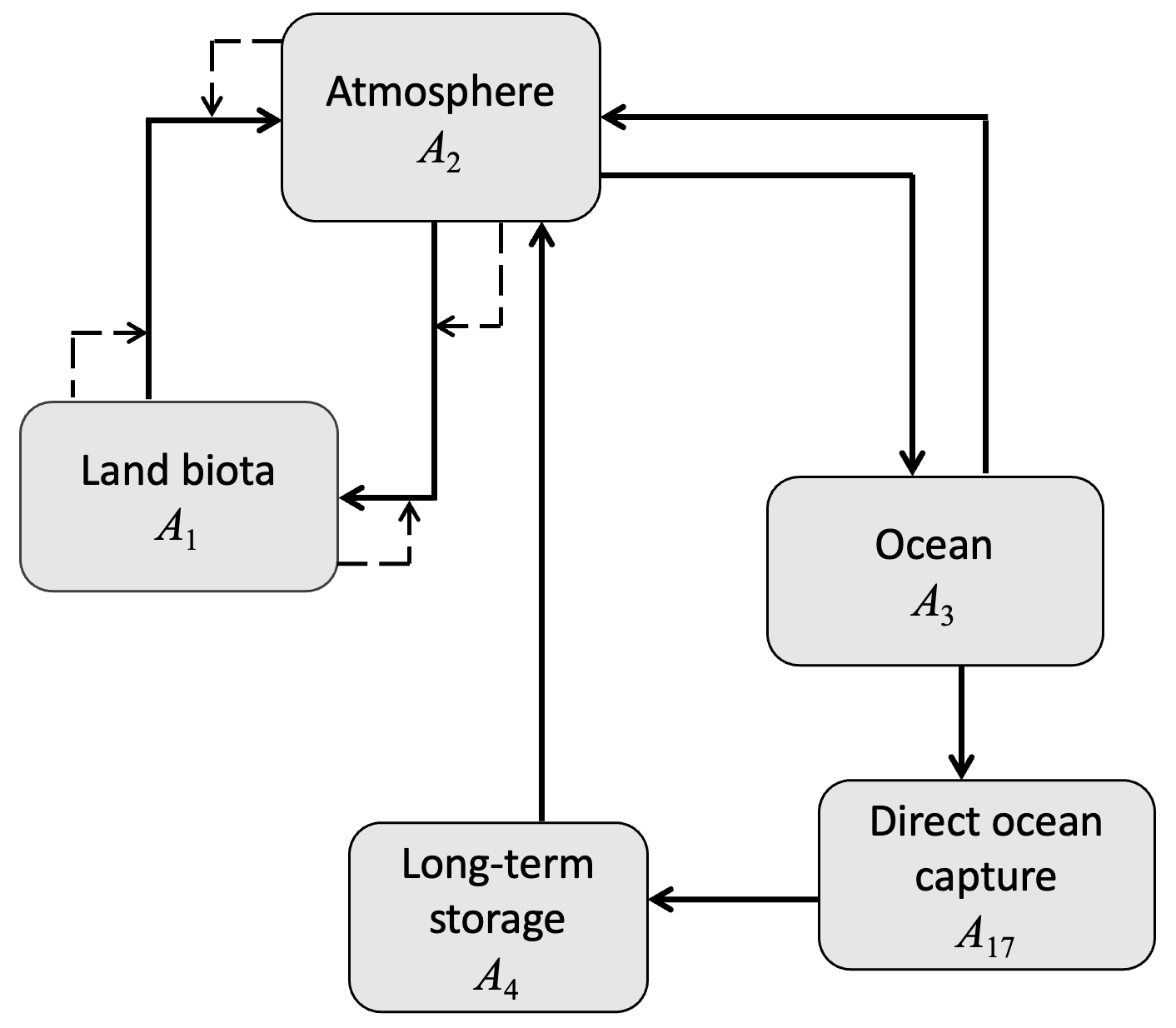}
    \caption{A biochemical map of the Earth's carbon cycle with direct ocean capture (DOC). The nodes represent the carbon pools. Furthermore, the solid arrows indicate carbon transfer, while the dashed arrows represent regulatory influences.}\label{fig:biochemicalmap}
    \end{center}
    \end{figure}

As seen in Figure \ref{fig:biochemicalmap},
the solid arrows indicate active carbon transfers between these pools, while the dashed arrows represent passive carbon transfers induced by 
regulatory influences.
For example, the solid arrow from the atmosphere ($A_2$) to the ocean ($A_3$) indicates that a portion of carbon in the atmosphere can be actively transferred to the ocean. Meanwhile, the transfer of carbon from land to atmosphere has both active and passive components and is influenced by both $A_1$ and $A_2$. As a result, we use dashed arrows in the diagram and write the reaction $A_1 + (A_1+A_2)\to A_2 + (A_1+A_2)$, which is the same as $2A_1+A_2 \to A_1+2A_2$, to represent both components of carbon transfer.

The modeling framework utilizes a power-law system in which the processes (or reactions) are represented by power-law functions. The structure of the rate function of the two processes with regulatory influences, where the two species $A_1$ and $A_2$ are involved, follows the form $k a_1^p a_2^q$ \cite{fortun2018deficiency}.

Our extended model includes an additional compartment ($A_4$) for the total carbon stock, which facilitates the transfer of carbon to the atmosphere at a linear rate. Furthermore, we incorporate an additional compartment for direct ocean capture ($A_{17}$), also with a linear rate. These extensions to the original Anderies model allow us to define four functional subsystems: the Anderies pre-industrial carbon cycle subsystem, the direct ocean capture subsystem, the carbon storage subsystem, and the carbon emission subsystem.

The biochemical map in Figure \ref{fig:biochemicalmap} can be represented as a chemical reaction network taking the different carbon pools as our species and the carbon transfers as reactions. The reactions in the DOC systems's corresponding network $\NN,$ together with the corresponding rate functions for each of the are given by
\begin{align*}
    &R_1: A_1 + 2A_2 \rightarrow 2A_1 + A_2 \quad \quad &&(k_1a_1^{p_1}a_2^{q_1}),\\
    &R_2: 2A_1 +A_2 \rightarrow A_1 + 2A_2 \quad \quad &&(k_2a_1^{p_2}a_2^{q_2}),\\
    &R_3: A_2 \rightarrow A_3 \quad \quad &&(k_3a_2), \\
    &R_4: A_3 \rightarrow A_2 \quad \quad &&(k_4a_3), \\
    &R_5: A_4 \rightarrow A_2 \quad \quad && (k_5a_4), \\
    &R_6: A_{17} \rightarrow A_4 \quad \quad &&(k_6a_{17)}, \\
    &R_7: A_3 \rightarrow A_{17} \quad \quad &&(k_7a_3),
\end{align*}
where the $k_i$'s, for $i = 1, 2, \ldots, 7,$ denote the rate constants for each of the seven reactions. We can also represent the underlying network of the DOC system in the following diagram

\[\begin{aligned}
    A_1 + 2A_2 \xrightleftharpoons{\qquad} A_2 + 2A_1\qquad \begin{tikzpicture}[node distance = {20mm}, baseline=(current  bounding  box.center)]
			\node (A17) {$A_{17}$};
            \node (A3) [above of=A17]{$A_3$};
            \node (A2) [left of=A3]{$A_2$};
            \node (A4) [left of=A17]{$A_4$};
            \draw[->] (A3) -- (A17);
            \draw[->] (A17) -- (A4);
            \draw[->] (A4) -- (A2);
            \draw [-left to] ($(A2.east) + (0pt, 2pt)$) -- ($(A3.west) + (0pt, 2pt)$);
    \draw [-left to] ($(A3.west) + (0pt, -2pt)$) -- ($(A2.east) + (0pt, -2pt)$);
        \end{tikzpicture}.
\end{aligned}\]

The kinetic order values of our system, as well as additional important quantities for our analysis are given in Table \ref{modelparameters}. In the table, we define the interaction differences of respiration and photosynthesis in the land biota and atmosphere, as well as their corresponding difference ratios, $R$ and $Q$.

\begin{table}[h!]
\centering
\resizebox{\textwidth}{!}{%
\begin{tabular}{|c|>{\raggedright\arraybackslash}p{10cm}|}
\hline
\textbf{Notation} & \textbf{Definition} \\ \hline
$p_1$ & kinetic order of land photosynthesis interaction ($p$-interaction) \\ \hline
$p_2$ & kinetic order of land respiration interaction ($r$-interaction) \\ \hline
$q_1$ & kinetic order of atmosphere photosynthesis interaction ($p$-interaction) \\ \hline
$q_2$ & kinetic order of atmosphere respiration interaction ($r$-interaction) \\ \hline
$p_2 - p_1$ & land $r$-$p$-interaction difference \\ \hline
$q_2 - q_1$ & atmosphere $r$-$p$-interaction difference \\ \hline
$R = \dfrac{p_2-p_1}{q_2-q_1}$ & land-atmosphere $r$-$p$-interaction difference ratio \\ \hline
$Q = \dfrac{q_2-q_1}{p_2-p_1}$ & atmosphere-land $r$-$p$-interaction difference ratio \\ \hline
\end{tabular}%
}
\caption{Model parameters in the DOC system}
\label{modelparameters}
\end{table}

From the defined parameters, we generate the ordinary differential equations (ODEs) that describe the dynamics of the network given by\begin{align*}
    \dfrac{da_1}{dt} &= k_1a_1^{p_1}a_2^{q_1} - k_2a_1^{p_2}a_2^{q_2} \\
    \dfrac{da_2}{dt} &= k_2a_1^{p_2}a_2^{q_2} - k_1a_1^{p_1}a_2^{q_1} - k_3a_2 + k_4a_3 + k_5 a_4  \\
    \dfrac{da_3}{dt} &= k_3a_2 - k_4a_3 - k_7 a_3  \\
    \dfrac{da_4}{dt} &= k_6a_{17} - k_5 a_4  \\
    \dfrac{da_{17}}{dt} &= k_7 a_3 - k_6a_{17}.
\end{align*}

We proceed by defining four different classes of the direct ocean capture system based on the signs of $R$ and $Q$. This concept of classifying carbon systems was introduced by Fortun and Mendoza in \cite{fortun2023comparative} and we provide a similar classification to our systems given the modification for direct ocean capture technology.

\begin{definition}
The set of direct ocean capture systems such that $R > 0$ ($R < 0$) is denoted by ${\sf DOC}_>$ (${\sf DOC}_<$). Elements of ${\sf DOC}_>$ (${\sf DOC}_<$) are called positive (negative) DOC systems.
\end{definition}

Notice that the ratios $R$ and $Q$ are multiplicative reciprocals of each other. Hence, both difference ratios must have the same signs. As such, we can equivalently check the sign of $Q$ to determine whether a DOC system is positive or negative. To be precise, a DOC system is also said to be positive (negative) if $Q > 0$ ($Q < 0$).

Now, note that the ratios $R$ and $Q$ are defined for $q_1 \neq q_2$ and $p_1\neq p_2,$ respectively. Moreover, $R$ and $Q$ are zero if $p_1 = p_2$ and $q_1 = q_2,$ respectively. This allows us to define two more classes for our DOC systems.

\begin{definition}
    The set of all direct ocean capture systems such that $R = 0,$ i.e. $p_1 = p_2$ but $q_1\neq q_2$ is denoted by ${\sf DOC}_{P_0},$ Similarly, if $Q = 0,$ i.e. $q_1 = q_2$ but $p_1\neq p_2,$ then this set is denoted by ${\sf DOC}_{Q_0}.$ Here, DOC systems in ${\sf DOC}_{P_0}$ (${\sf DOC}_{Q_0}$) are said to be a $P$-null ($Q$-null) DOC system.
\end{definition}

As much as possible, we formulate our results in this study in terms of the four classes of DOC systems, namely, ${\sf DOC}_>, {\sf DOC}_<, {\sf DOC}_{P_0},$ and ${\sf DOC}_{Q_0}.$ The following sections discuss the existence, multiplicity, and absolute concentration robustness (ACR) in the four classes of DOC systems and, if necessary, some of its specific subsets. 

{
\subsection{Existence of positive steady states in the DOC model}
\label{subsec:existence}

Steady states typically describe the long-term behaviors of chemical or biochemical systems. Mathematically, at these states, the time derivatives vanish, meaning that over a long period of time, the concentrations of the species remain constant.

As we study the steady state related properties of our DOC model, we first check the existence of its positive steady states depending on the rate constants. Since the model follows power law kinetics, we can use the results presented by Alamin and Hernandez in \cite{alamin2024positive} to verify the existence of positive steady states for the entire system through the subsystems induced by its underlying independent subnetworks (see Section \ref{subs:network:decomposition} for details). Applying the MATLAB program \cite{LubeniaINDECS} to our model for finding independent decompositions gives two independent subnetworks ($\mathcal N_1$ and $\mathcal N_2$) whose reactions sets are given by $\mathcal R_1 = \{R_1, R_2\}$ and $\mathcal R_2 = \{R_3, R_4, R_5, R_6, R_7\},$ respectively. 

Regardless of the interactions of photosynthesis and respiration on land and in the atmosphere, the decomposition of the network into its connected components satisfies both requirements in checking the existence of positive steady states: stoichiometric independence and $\widehat T$-independence (see Appendix \ref{details:existence} for details). Therefore, invoking Theorem 2 in \cite{alamin2024positive}, all defined classes of our DOC systems, i.e. positive, negative, $P$-null, and $Q$-null, have a positive steady state for any rate constants if and only if each subsystem induced by the independent subnetworks of the decomposition also has a positive steady state for any set of rate constants. 

Indeed, whenever $p_1 \neq p_2$ or $q_1 \neq q_2,$ we can use the method of Johnston et al. \cite{johnston2019deficiency}, the steps of \cite{hernandez2023framework, hernandez2023network}, as well as the computational package COMPILES from \cite{hernandez2023framework} to show that the set of positive steady states of the independent subsystems of the DOC system is nonempty and may be parametrized in terms of its rate constants (see Appendix \ref{details:parametrization} for details). In terms of our defined classes, this means that all classes of our DOC systems admit at least one positive steady state for any set of rate constants. On the other hand, when both $p_1 = p_2$ and $q_1= q_2$, 
the corresponding ODEs of the first subnetwork are
\begin{align*}
    \dfrac{da_1}{dt} &= k_1a_1^{p_1}a_2^{q_1} - k_2a_1^{p_2}a_2^{q_2} = k_1a_1^{p_1}a_2^{q_1} - k_2a_1^{p_1}a_2^{q_1}=(k_1-k_2)a_1^{p_1}a_2^{q_1},\\
    \dfrac{da_2}{dt} &= k_2a_1^{p_2}a_2^{q_2} - k_1a_1^{p_1}a_2^{q_1}=k_2a_1^{p_1}a_2^{q_1} - k_1a_1^{p_1}a_2^{q_1}=(k_2-k_1)a_1^{p_1}a_2^{q_1}.
\end{align*}
Thus, positive steady states of this subnetwork's associated system exist only when $k_1 = k_2.$ For independent subnetworks, the intersection of their sets of positive steady states is the set of positive steady states of the whole network (see Theorem \ref{thm:feinberg-decomp} for details). Since the decomposition is also $\widehat T$-independent, we have the existence of positive steady states for the entire system only if $k_1 = k_2$ \cite{alamin2024positive}. Furthermore, the DOC systems admit a positive steady state for every stoichiometric class, i.e., for each positive value of the conserved total quantity (see Appendix \ref{details:stoichiometric:class} for details).

\subsection{Conditions for multistationarity of the DOC system}
\label{subsec:multi}

The capacity of the DOC system to admit multiple steady states depends on the values of its kinetic orders. 
 Specifically, whether the system can exhibit multistationarity is determined by the sign of its interaction difference ratios $$R=\dfrac{p_2-p_1}{q_2-q_1} \text{ and } Q=\dfrac{q_2-q_1}{p_2-p_1}.$$

For the DOC system, the stoichiometric subspace $S$ is given as follows (see Apppendix \ref{details:multistationarity}  for details):
$$S=\rm{span} \left \{\begin{bmatrix} 1 \\ -1 \\ 0 \\ 0 \\ 0 \end{bmatrix}, \begin{bmatrix} 0 \\ -1 \\ 1 \\ 0 \\ 0 \end{bmatrix}, \begin{bmatrix} 0 \\ 0 \\ -1 \\ 0 \\ 1 \end{bmatrix}, \begin{bmatrix} 0 \\ 0 \\ 0 \\ 1 \\ -1 \end{bmatrix} \right \}.$$
Furthermore, the orthogonal complement $(\tilde{S})^{\perp}$ of the kinetic flux subspace $\tilde{S}$ is given as follows (see Apppendix \ref{details:multistationarity} for details):
\[ (\tilde{S})^{\perp} =
\begin{cases*}
 \rm{span} \left \{\begin{bmatrix} -Q \\ 1 \\ 1 \\ 1 \\ 1 \end{bmatrix} \right \} & when written in terms of $Q$\\
\rm{span} \left \{\begin{bmatrix} -1 \\ R \\ R \\ R \\ R \end{bmatrix} \right \} & when written in terms of $R$  
\end{cases*}.
\]
The subspaces $S$ and $(\tilde{S})^{\perp}$ are needed in the simple criterion of M\"uller and Regensburger \cite{regensburger}, which determines when a system admits more than one (complex balanced) steady state. This is performed by examining the possible sign patterns of the vectors in these spaces, which in turn rely on the mentioned difference ratios. Using such criterion (see Appendix \ref{details:multistationarity} for details), the system has the capacity to admit multiple steady states when $R > 0.$ Therefore, all positive DOC systems are multistationary and can admit multiple steady states under the same set of parameters.

Next, we use the criterion of Feliu and Wiuf \cite{interacting,wiuf} and individually investigate the signs of kinetic orders $p_1,p_2,q_1$ and $q_2$ to conclude when the system is monostationary (see Apppendix \ref{details:multistationarity} for details). Indeed, by this criterion (see Theorem \ref{powerlawinjectivity}), the DOC system admits a unique positive steady state whenever $p_1, q_2 < 0$ and $p_2, q_1 > 0$, or whenever $p_1, q_2 > 0$ and $p_2, q_1 < 0$. Note that the conditions on the signs of the kinetic orders imply that $R < 0$ but such negativity condition for the difference ratios themselves are not sufficient conditions for the DOC system to become monostationary. We note then that although all monostationary DOC systems are negative, the converse does not immediately follow. Only systems in the subset of ${\sf DOC}_{<}$ satisfying $p_1,q_2 < 0$ and $p_2, q_1 > 0$ or $p_1, q_2 > 0$ and $p_2, q_1 < 0$ have been shown to exhibit monostationarity.

\subsection{Conditions for ACR of the DOC system}
\label{subsec:acr}

Following the method of Hernandez et al. for computing steady state parametrization of chemical reaction networks \cite{hernandez2023network,hernandez2023framework}, we obtain some positive steady state parametrizations of the DOC system for different values for the kinetic orders (see Appendix \ref{details:parametrization} for the computation). First, when the kinetic order of the land photosynthesis interaction differs from that of the land respiration interaction, that is, when $p_1 \ne p_2,$ we have
\begin{align*}
a_1 &= \left(\dfrac{k_1}{k_2}\right)^{\frac{1}{p_2 - p_1}} \tau^{q_2 - q_1} \\
a_ 2 &= \tau^{p_1 - p_2}\\
a_3 &= \dfrac{k_3}{k_4 + {k_7}} \tau^{p_1 - p_2} \\
a_4 &= \dfrac{k_3 {k_7}}{{k_5}(k_4 + {k_7})} \tau^{p_1 - p_2} \\
a_{17} &= \dfrac{k_3 {k_7}}{k_6 (k_4 + {k_7})} \tau^{p_1 - p_2}.
\end{align*}
where $\tau>0$. To ensure that the concentration of $\text{CO}_2$ in land ($A_1$) remains stable regardless of the initial conditions, the kinetic order of the atmosphere photosynthesis interaction should be equal to the kinetic order of the atmosphere respiration interaction, i.e., $q_1 = q_2$. The concentration of carbon dioxide in the other carbon pools remains variable in this scenario. Note that in this case, the atmosphere-land $r$-$p$-interaction difference ratio becomes zero, i.e., $Q = 0.$ Thus, we achieve absolute concentration robustness (ACR) on $A_1$ only whenever we have a $Q$-null DOC system.

Analogously, when the kinetic order of the atmosphere photosynthesis interaction is different from the kinetic order of the atmosphere respiration interaction, i.e., $q_1 \ne q_2$, we obtain the following parametrization of positive steady states of the DOC system:
\begin{align*}
    a_1 &= \tau^{q_2-q_1} \\
    a_2 &= \left(\dfrac{k_1}{k_2}\right)^{\frac{1}{q_2-q_1}} \tau^{p_1-p_2} \\
    a_3 &= \dfrac{k_3}{k_4 + {k_7}} \left(\dfrac{k_1}{k_2}\right)^{\frac{1}{q_2-q_1}} \tau^{p_1-p_2} \\
    a_4 &= \dfrac{k_3{k_7}}{{k_5}(k_4+{k_7})} \left(\dfrac{k_1}{k_2}\right)^{\frac{1}{q_2-q_1}} \tau^{p_1-p_2} \\
    a_{17} &= \dfrac{k_3{k_7}}{k_6(k_4+{k_7})} \left(\dfrac{k_1}{k_2}\right)^{\frac{1}{q_2-q_1}} \tau^{p_1-p_2}
\end{align*}
where $\tau>0$. Similar to our analysis in the previous case, the concentrations of CO$_2$ in the atmosphere ($A_2$), ocean ($A_3$), DOC ($A_{17}$), and total carbon stock pools ($A_4$) is stable whenever the kinetic order of the land photosynthesis interaction equals that of the land respiration interaction, i.e. $p_1 = p_2.$ In this case, we have $R = 0$ and so we achieve ACR in species $A_2, A_3, A_4,$ and $A_{17}$ in $P$-null DOC systems.

Finally, if we have $p_1 \neq p_2$ and $q_1\neq q_2,$ then we can take any of the above positive steady state parametrizations for our system. As a result, the concentration of all species at their respective steady states will vary over different sets of initial concentrations. Therefore, for the case of positive and negative DOC systems, we do not achieve ACR in any of the species.

\begin{remark}
    An alternative method to determine the ACR property of the DOC system is using the species hyperplane criterion \cite{hyperplane:criterion}. This states that a system has ACR species if and only if the vector coordinates corresponding to these species are zero for all basis vectors in $(\tilde{S})^{\perp}$. Recall from the previous section that 
    \[ (\tilde{S})^{\perp} =
\begin{cases*}
 \rm{span} \left \{\begin{bmatrix} {-Q} \\ 1 \\ 1 \\ 1 \\ 1 \end{bmatrix} \right \} & when written in terms of ${Q}$\\
\rm{span} \left \{\begin{bmatrix} -1 \\ {R} \\ {R} \\ {R} \\ {R} \end{bmatrix} \right \} & when written in terms of {$R$}  
\end{cases*}.
\]
Hence, positive and negative DOC systems has no ACR in any species. On the other hand, $P$-null systems, with $R=0$, have ACR in species $A_2$, $A_3$, $A_4$ and $A_{17}$, while $Q$-null systems, with $Q=0$, have ACR in species $A_1$.
\end{remark}


\subsection{Sufficient conditions for carbon reduction}
\label{subsec:suff-cond}

In this section, we determine sufficient conditions to ensure that for any set of initial conditions $A_1^{\circ}$, $\ldots$, $A_4^{\circ}$, $A_{17}^{\circ}$ and any set of steady state values $A_1^{*}$, $\ldots$, $A_4^{*}$, $A_{17}^{*}$ in an associated stoichiometric class $S^{\circ}$, there is carbon pool reduction in the ocean, i.e., $A_3^{\circ}>A_3^{*}$.
In other words, the long-term concentration of carbon in the ocean is lower than its initial concentration.
The approach taken here is to use the conserved quantity of the DOC, i.e., $$T=A_1^{\circ}+A_2^{\circ}+A_3^{\circ}+A_4^{\circ}+A_{17}^{\circ}.$$

Note that the underlying network of the DOC system is conservative, and hence each stoichiometric class is compact \cite{gmak}. Hence, we can define the continuous maps pr$_{i}$: $\mathbb{R}^m \to \mathbb{R}$ where $i$ denotes the index of the carbon pools of our system, i.e. $\mathrm{pr}_2(A)$ is the concentration of $A_2$ in the system. Note that these maps and their sums attain maxima and minima on any of its stoichiometric class or closed subset. We now present a sufficient condition on the network parameters for carbon reduction in the ocean.

\begin{proposition}
    Suppose a DOC system has initial conditions $A_i^\circ$ and steady state values $A_i^*$ in the associated stoichiometric class $S^\circ.$ Let $m'$ be the minimum of ${\rm pr}_2$ and $M'$ be the maximum of ${\rm pr}_1 + {\rm pr}_2 + {\rm pr}_4 + {\rm pr}_{17}$ on $S^\circ.$ Then $A_3^* < A_3^\circ$ whenever $\dfrac{k_3}{k_4 + k_7} < \dfrac{T - M'}{m'}.$

    \begin{proof}
        We consider the cases when $p_1 \neq p_2$ or $q_1 \neq q_2.$ These two cases are sufficient to describe the behavior for positive, negative, $P$-null, and $Q$-null systems. Note that if either $p_1 = p_2$ or $q_1 = q_2,$ but not both, steady states of some species may be parametrized by $A_2$. Specifically, for systems where $p_1\neq p_2,$ i.e. $Q$-null if $q_1 = q_2$ and positive/negative otherwise, we have \begin{align*}
            A_1 &= \left(\dfrac{k_1}{k_2}\right)^{\frac{1}{p_2 - p_1}} \\
A_2 &= A_2\\
A_3 &= \dfrac{k_3}{k_4 + k_7}A_2\\
A_4 &= \dfrac{k_3 k_7}{k_5(k_4 + k_7)} A_2 \\
A_{17} &= \dfrac{k_3 k_7}{k_6 (k_4 + k_7)}A_2.
        \end{align*} For systems with $q_1\neq q_2,$ i.e. $P$-null if $p_1 = p_2$ and positive/negative otherwise, we have \begin{align*}
    A_1 &= \tau^{q_2-q_1} \\
    A_2 &= A_2 \\
    A_3 &= \dfrac{k_3}{k_4 + {k_7}} A_2 \\
    A_4 &= \dfrac{k_3 {k_7}}{{k_5}(k_4+{k_7})} A_2 \\
    A_{17} &= \dfrac{k_3{k_7}}{k_6(k_4+{k_7})} A_2.
\end{align*} For either case, we have $A_3 = \dfrac{k_3}{k_4 + k_7}A_2.$ Thus, we have \[A_3^* = \frac{k_3}{k_4 + k_7}A_2^* < \frac{T - M'}{m'}A_2^* \leq T - M' \leq T - (A_1^\circ + A_2^\circ + A_4^\circ + A_{17}^\circ) = A_3^\circ,\] which gives our desired result.
    \end{proof}
\end{proposition}

We compare the sufficient conditions for carbon reduction for our proposed DOC model to one with direct air capture. In contrast to our DOC model which looks at the sufficient conditions for carbon reduction in the ocean, we look at the sufficient conditions for carbon reduction in the atmosphere for models with direct air capture. We present these sufficient conditions in Table \ref{tab:compare-carbon-reduction-suff}. We remark that the sufficient conditions for carbon reduction in the atmosphere with DAC involve more parameters than those of carbon reduction in the ocean with DOC.

\subsection{Tabular summary of dynamic properties of the DOC system}
In this work, we study the long-term behaviors of a global carbon cycle model incorporating Direct Ocean Capture (DOC) technology through its positive steady states. Applying results in chemical reaction network theory, we were able to provide conditions for the existence of positive steady states, multistationarity, and absolute concentration robustness in our DOC model.

First, we have shown that for all four defined classes of DOC systems, namely, the positive, negative, $P$-null, and $Q$-null systems, there exists at least one positive steady state for any set of rate constants. In contrast, if a DOC system is not in any of these classes, i.e. when both $p_1 = p_2$ and $q_1 = q_2,$ positive steady states exist only when $k_1 = k_2.$

\begin{landscape}
\begin{table}
    \centering
    \scriptsize
    \renewcommand{\arraystretch}{1.5}
    \begin{tabular}{lcc}
        \hline
        & Carbon reduction in ocean & Carbon reduction in atmosphere \\
        & ($A_3^* < A_3^\circ$) & ($A_2^* < A_2^\circ$) \\
        \hline
        Positive ($R > 0$) 
        & \multirow{4}{*}{$\dfrac{k_3}{k_4 + k_7} < \dfrac{T - M'}{m'}$} 
        & \multirow{2}{*}{$1 + \dfrac{M''}{m'} < \left(\dfrac{k_1}{k_2}\right)^{\frac{1}{p_2 - p_1}}(m')^{-Q} + k^\star$ or $1 + \dfrac{M''}{m'} < \left(\dfrac{k_1}{k_2}\right)^{\frac{1}{q_2 - q_1}}(m')^{-R} + k^\star$} \\
        Negative ($R < 0$) & & \\
        $P$-null ($P = 0$) & & \multirow{2}{*}{$1 + \dfrac{M''}{m'} < k^\star$} \\
        $Q$-null ($Q = 0$) & & \\
        \hline
        \multicolumn{3}{l}{where $k^\star = \dfrac{k_3k_5k_7 + k_6k_4(k_5 + k_7)}{k_4k_5k_7}$} \\
    \end{tabular}
    \caption{Sufficient conditions for carbon reduction in carbon cycles with DOC and DOC technologies.}
    \label{tab:compare-carbon-reduction-suff}
\end{table}
\end{landscape}

Next, we have shown that all positive DOC systems can admit more than one positive steady states for a fixed set of parameters, i.e. they are multistationary. In contrast, not all negative DOC systems exhibit monostationarity. Specifically, only {{two subsets}} of our negative DOC systems, satisfying $p_1, q_2 < 0$ and $p_2, q_1 > 0$, {{or $p_1, q_2 > 0$ and $p_2, q_1 < 0$}} admit a unique positive steady state for each fixed set of parameters. We emphasize that although monostationarity in a DOC system implies that the system is negative, the converse does not follow. For $P$-null and $Q$-null DOC systems, we investigate the induced ODEs and associated conservation laws of these systems to determine the multiplicity of their steady states. Following this approach, we conclude that all ${\sf DOC}_{P_0}$ and ${\sf DOC}_{Q_0}$ systems are monostationary (see Theorems \ref{thm-monoP0} and \ref{thm-monoQ0} for the proofs). These results are validated in Appendix \ref{appendix:validation:multi}.

Finally, absolute concentration robustness, or ACR, on some species of our DOC system was exhibited for $P$-null and $Q$-null systems only. Specifically, regardless of initial concentrations, stable concentrations at steady state for the atmosphere ($A_2$), ocean ($A_3$), direct ocean capture ($A_{17}$), and total carbon stock pools ($A_4$), are achieved in $P$-null DOC systems. For $Q$-null DOC systems, this allows for ACR of carbon dioxide in the land biota ($A_1$) only. These results are validated in Appendix \ref{appendix:validation:ACR}.

\begin{table}
\centering
\resizebox{\textwidth}{!}{%
\begin{tabular}{p{4cm} p{10cm}}
    \hline
    Steady State Property & Class of DOC systems* \\
    \hline
    Existence
    & ${\sf DOC}_>$\,\,: for any set of rate constants \\
    &
    ${\sf DOC}_<$\,\,: for any set of rate constants \\
    &
    ${\sf DOC}_{P_0}$\,: for any set of rate constants \\
    &
    ${\sf DOC}_{Q_0}$: for any set of rate constants \\
    
    Multiplicity
    & ${\sf DOC}_>$\,\,: multistationary \\
    &
    ${\sf DOC}_<$\,\,: contains monostationary systems; \\
    &
    \qquad\qquad must satisfy $p_1, q_2 < 0$ and $p_2, q_1 > 0$ \\
    &
    \qquad\qquad or $p_1, q_2 > 0$ and $p_2, q_1 < 0$ \\
    &
    ${\sf DOC}_{P_0}$\,: monostationary \\
    &
    ${\sf DOC}_{Q_0}$: monostationary \\
    
    ACR
    & ${\sf DOC}_>$\,\,: no ACR \\
    &
    ${\sf DOC}_<$\,\,: no ACR \\
    &
    ${\sf DOC}_{P_0}$: ACR in $A_2, A_3, A_4,$ and $A_{17}$ only \\
    &
    ${\sf DOC}_{Q_0}$: ACR in $A_1$ only \\
    \hline
\end{tabular}
}
\vspace{0.5em}

\footnotesize
* Positive (${\sf DOC}_>$), Negative (${\sf DOC}_<$), $P$-null (${\sf DOC}_{P_0}$), or $Q$-null (${\sf DOC}_{Q_0}$).

\caption{Summary of steady state properties of the different classes of DOC systems}
\label{tab:summary}
\end{table}

Ideally, carbon dioxide removal through direct ocean capture technology aims to achieve a stable and unique concentration of carbon dioxide in the ocean at its steady state. Our results imply that this ideal situation can be achieved in $P$-null direct ocean capture systems.

\subsection{Comparison of carbon capture systems}

\begin{figure}[H]
    \begin{center}
    \includegraphics[width=11cm,height=7cm,keepaspectratio]{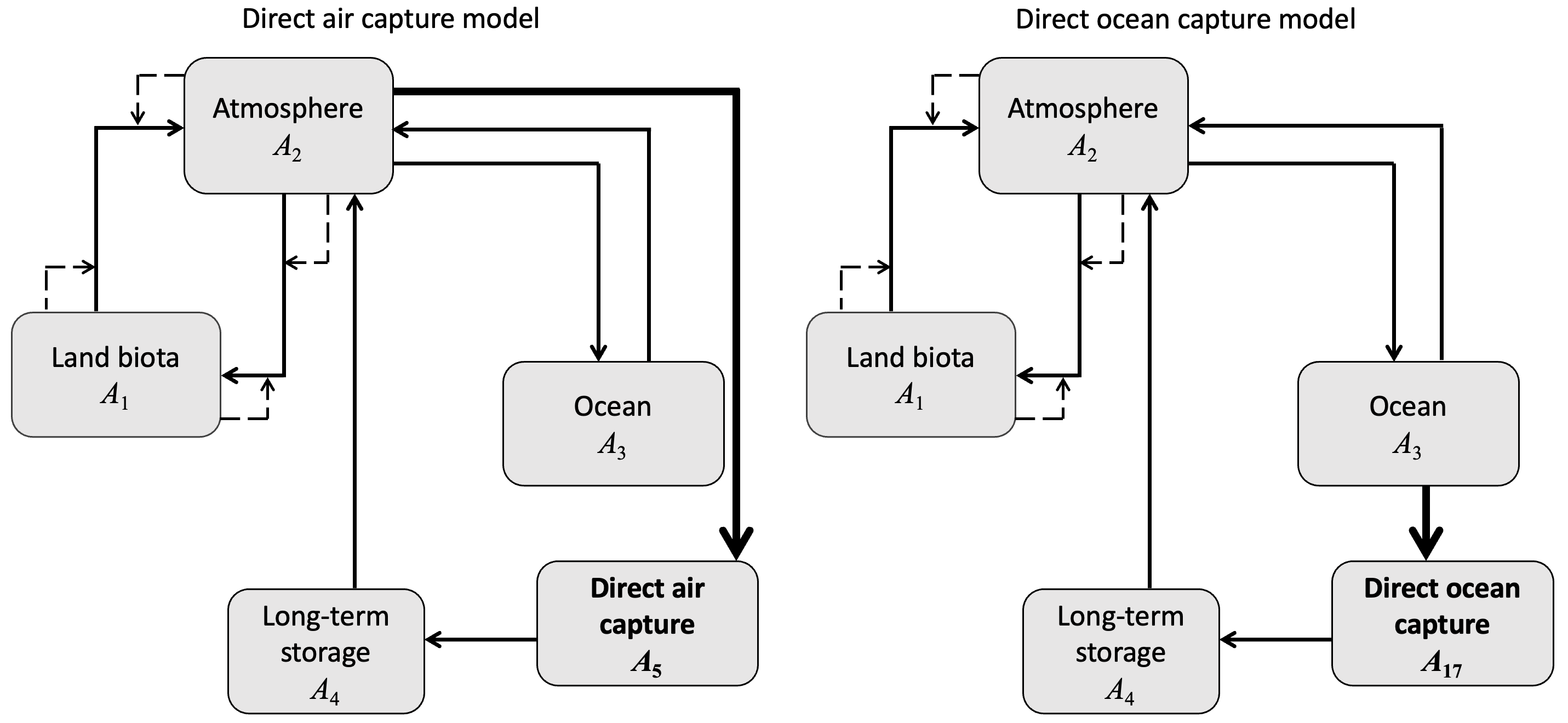}
    \caption{{Side-by-side comparison of the biochemical maps of the underlying networks of DAC and DOC systems, placed on the left and right panels, respectively. The differences between the two networks are highlighted in red.}}\label{fig:DAC:DOC}
    \end{center}
    \end{figure}

In this section, we begin to compare carbon capture systems, focusing specifically on models which utilize either direct air capture (DAC) or direct ocean capture (DOC) exclusively. Figure \ref{fig:DAC:DOC} shows the biochemical maps of the underlying networks of the DAC and DOC systems. Both models contain the four basic carbon pools ($A_1, A_2, A_3$ and $A_4$), along with their associated processes. The key difference of the models is the exclusive presence of compartments $A_5$ (direct air capture) and $A_{17}$ (direct ocean capture) in the DAC and DOC models, respectively. Furthermore, the reaction that captures the carbon from the atmosphere ($A_2 \to A_5$) is present only in the DAC network, whereas the reaction that captures the carbon from the ocean ($A_3 \to A_{17}$) is present only in the DOC network. In Table \ref{tab:NetworkNumbers}, we list and compare some network numbers describing network structure and composition to compare the DAC and DOC models. Using the standard CRNToolbox \cite{FeinbergToolbox} to obtain these numbers, we note that the network numbers for both models match exactly. 

\begin{table}[h!]
    \centering
    \begin{tabular}{|l|c|c|c|}
    \hline
    \textbf{Network Numbers} & \textbf{Notation} & \textbf{DAC} & \textbf{DOC}\\
    \hline
    Species & $m$ & $5$ & $5$ \\
    \hline
    Complexes & $n$ & $6$ & $6$ \\
    \hline
    Reactant Complexes & $n_r$ & $6$ & $6$\\
    \hline
    Reversible Reactions & $r_{\text{rev}}$ & $2$ & $2$ \\
    \hline
    Irreversible Reactions & $r_{\text{irrev}}$ & $3$ & $3$ \\
    \hline
    Reactions & $r$ & $7$ & $7$ \\
    \hline
    Linkage Classes & $\ell$ & $2$ & $2$ \\
    \hline
    Strong Linkage Classes & $s\ell$ & $2$ & $2$  \\
    \hline
    Terminal Strong Linkage Classes & $t$ & $2$ & $2$ \\
    \hline
    Rank & $s$ & $4$ & $4$ \\
    \hline
    Deficiency & $\delta$ & $0$ & $0$ \\
    \hline 
    \end{tabular}
    \caption{Network numbers of the models exclusively with DAC and DOC technologies using CRNToolbox}
    \label{tab:NetworkNumbers}
\end{table}
We also obtain from CRNToolbox \cite{FeinbergToolbox} the coincidence of the some structural properties of the DOC and DAC models as seen in Table \ref{tab:CRNToolbox-DAC-DOC}.
\begin{table}[ht!]
\centering
\resizebox{\textwidth}{!}{%
\begin{tabular}{|l|c|c|l|}
    \hline
    \multicolumn{1}{|l|}{Property} & DAC & DOC & Description of the property \\
    \hline
    Deficiency zero & Yes & Yes & The deficiency is a non-negative integer that \\
                   &     &     & measures the linear dependence of the reactions. \\
    \hline
    Weakly reversible & Yes & Yes & Each reaction belongs to a cycle. \\
    \hline
    Positive dependent & Yes & Yes & There is a set of positive numbers for which the \\
                       &     &     & linear combination of the reaction vectors in the \\
                       &     &     & network equals zero. \\
    \hline
    Conservative & Yes & Yes & There is a vector in the positive orthant that is \\
                &     &     & orthogonal to all the reaction vectors, hence, \\
                &     &     & respecting a conservation law. \\
    \hline
    Concordant & No & No & A structural property that enforces a degree of \\
              &    &    & dull, reliable behavior even against varieties \\
              &    &    & of kinetics; multistationarity is not possible. \\
    \hline
    Independent linkage classes & Yes & Yes & The linkage class decomposition is independent. \\
    \hline
    Maximally closed & Yes & Yes & The dimension of the stoichiometric subspace \\
                    &     &     & is one less than the number of species, i.e., $s = m - 1$. \\
    \hline
    High reactant diversity & Yes & Yes & The number of reactant complexes is more \\
                            &     &     & than the dimension of the stoichiometric \\
                            &     &     & subspace, i.e., $n_r > s$. \\
    \hline
\end{tabular}
}
\caption{Structural properties of the DAC and DOC networks obtained from the standard CRNToolbox}
\label{tab:CRNToolbox-DAC-DOC}
\end{table}

Finally, using the results of Fortun et al. \cite{fortun2024determining} and Table \ref{tab:summary} in this paper, we see a coincidence in the dynamic properties of the DAC-only and DOC-only models. Indeed, for all four classes that we have previously identified (i.e., positive, negative, $P$-null, and $Q$-null), both systems incorporating DAC and DOC exclusively exhibit the same dynamic properties on the existence of positive equilibrium, multistationarity and ACR.

Finally, in Section \ref{subsec:suff-cond}, we presented the sufficient conditions for oceanic carbon reduction for systems with DOC and atmospheric carbon reduction for systems with DAC. Notably, the sufficient conditions for the DAC model are more complex than that for the DOC model. This increased complexity may be attributed to the difference of the number of reactions occurring in the carbon compartments for the atmosphere and the ocean. For the DOC model, our results is consistent with our intuition that a higher value for the rate constants corresponding to an outflow of carbon in the ocean pool lead to greater oceanic carbon reduction efficiency. Similarly, for systems with direct air capture, rate constants corresponding to an outflow of carbon in the atmopshere also implies greater atmospheric carbon reduction efficiency.

\subsection{Analysis of the integrated air and ocean carbon capture system}

In this last section, we present a model which integrates both direct air capture and direct ocean capture technologies. In contrast to our comparison of the DAC-only and DOC-only systems, we show here that there are some differences to the network numbers and dynamic properties of the integrated carbon capture model.

First, we present the biochemical map of the integrated system with both DAC and DOC technologies in Figure \ref{fig:combined}. We note that this model now has six species, as a result of the integration of the carbon capture technologies.

\begin{figure}[H]
    \begin{center}
    \includegraphics[width=12cm,height=7cm,keepaspectratio]{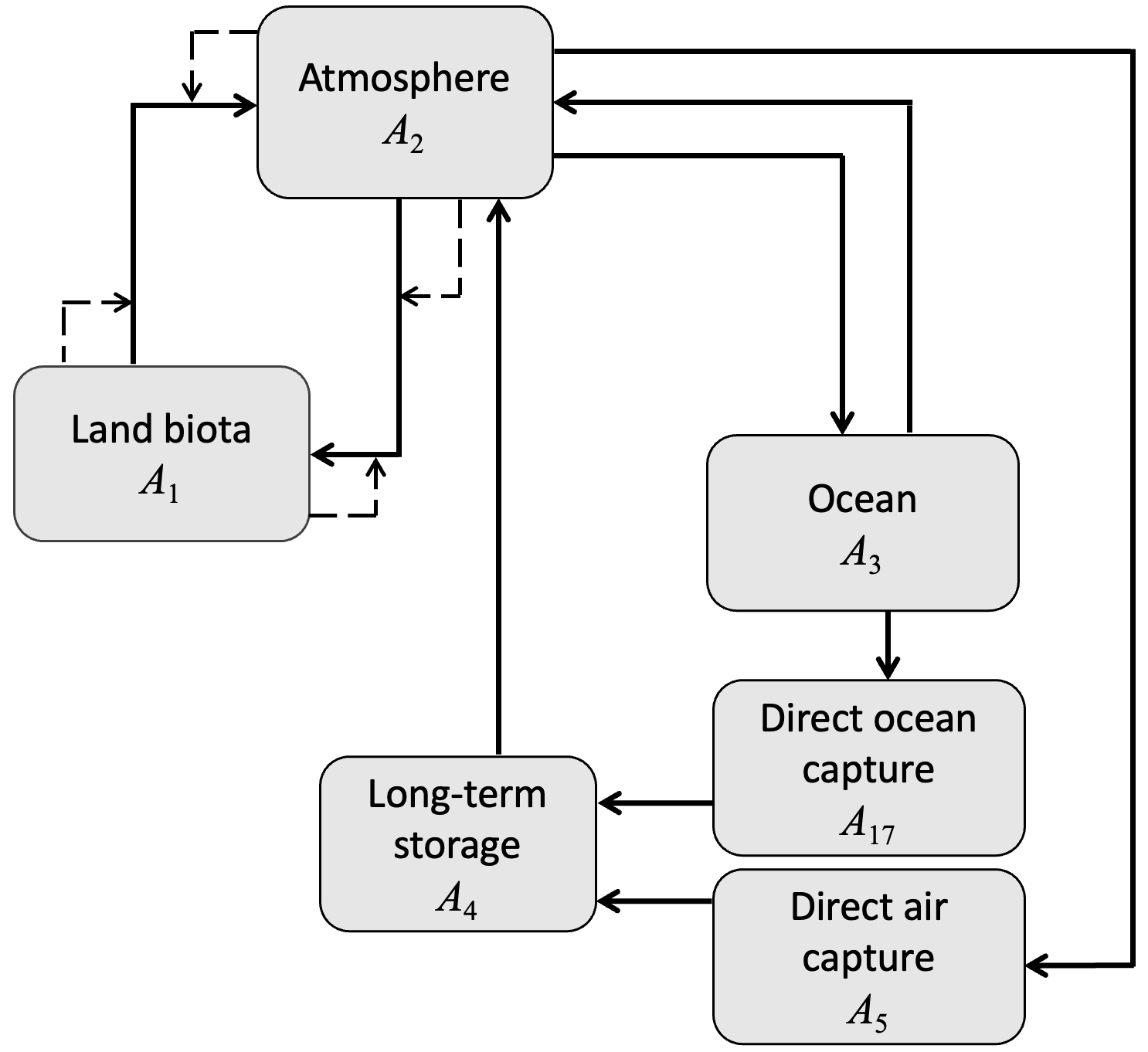}
    \caption{Biochemical map of the underlying networks of the integrated DAC and DOC systems.}\label{fig:combined}
    \end{center}
    \end{figure}

Next, Table \ref{tab:NetworkNumbers-Combined} gives the network numbers of the integrated system. The differences in some network numbers are immediate, but we do note that similar to the systems employing only one of the carbon capture technologies, the underlying network of the integrated system also has zero deficiency. In contrast, the rank of the integrated system is greater than one compared to the DAC-only and DOC-only due to the addition of the carbon pools for DOC and DAC, respectively.

\begin{table}[h!]
    \centering
    \begin{tabular}{|l|c|c|}
    \hline
    \textbf{Network Numbers} & \textbf{Notation} & \textbf{Integrated} \\
    \hline
    Species & $m$ & $6$ \\
    \hline
    Complexes & $n$ & $7$ \\
    \hline
    Reactant Complexes & $n_r$ & $7$ \\
    \hline
    Reversible Reactions & $r_{\text{rev}}$ & $2$ \\
    \hline
    Irreversible Reactions & $r_{\text{irrev}}$ & $5$ \\
    \hline
    Reactions & $r$ & $9$ \\
    \hline
    Linkage Classes & $\ell$ & $2$ \\
    \hline
    Strong Linkage Classes & $s\ell$ & $2$ \\
    \hline
    Terminal Strong Linkage Classes & $t$ & $2$ \\
    \hline
    Rank & $s$ & $5$ \\
    \hline
    Deficiency & $\delta$ & $0$ \\
    \hline 
    \end{tabular}
    \caption{Network numbers of the integrated network using CRNToolbox}
    \label{tab:NetworkNumbers-Combined}
\end{table}

The structural properties of the integrated system are the same with those of the DOC-only and DAC-only models. That is, the integrated carbon capture system satisfies the descriptions of deficiency zero, weakly reversible, positive dependent, conservative, independent linkage classes, maximally closed, and high reactant diversity, and is also not concordant. Table \ref{tab:CRNToolbox-DAC-DOC} shows the precise descriptions of these properties.

Finally, we now compare the dynamical properties of the integrated system with the DOC-only and DAC-only systems. Specifically, we compare the dynamical properties of the integrated system with the DOC-only system. For existence, the addition of the species for direct air capture does not affect the conditions required for a positive steady state. Instead, this addition increases the system's rank by one, which still allows it to satisfy the criteria established in \cite{alamin2024positive} for the existence of a positive steady state through its independent subnetworks. Moreover, the conditions for multiplicity of the integrated system is also the same as that of the DOC-only system, as shown in Appendix \ref{details:multi-combi}. Finally, based on the parametrization of the positive steady states of the integrated system (see Appendix \ref{details:param-combi}), the addition of the new species in the integrated system does not eliminate any species exhibiting ACR in any of the identified classes. Essentially, the set of species exhibiting ACR across four classes in the integrated system is the union of the same species exhibiting ACR in the DOC-only and DAC-only systems.

\section{Conclusion and Recommendation}
In this work, we study a global carbon cycle model that incorporates direct ocean capture (DOC) technology, using tools and concepts from Chemical Reaction Network Theory (CRNT) to analyze the system’s long-term behavior without requiring specific parameter values. Specifically, we investigate the existence and multiplicity of steady states and identify parameter conditions under which long-term concentration robustness emerges among the carbon pools.

Our analysis shows that the DOC model always admits at least one positive steady state for any set of rate constants. This implies the ability of the system to go to a nonzero concentration of carbon across its pools in the long term. We also derive conditions on the rate constants that give rise to multistationarity. In the context of the global carbon cycle, the presence of such multiple positive steady states may correspond to low- or high-carbon equilibria. A low-carbon steady state typically reflects a stable climate that supports biodiversity, agriculture, and habitability. In contrast, a high-carbon steady state may lead to global warming, extreme weather events, and ecosystem disruption.

Additionally, we identify conditions under which the system exhibits absolute concentration robustness. This study suggests unique and robust oceanic carbon concentrations arise when the kinetic order of atmospheric photosynthesis differs from that of atmospheric respiration.

We further extend our analysis by introducing a model that integrates both direct air capture (DAC) and direct ocean capture (DOC) technologies.
We find that this integration does not affect the required conditions for on positive steady state.
In the combined system, the set of species exhibiting ACR across four different model classes corresponds to the union of the ACR species identified in the DAC-only and DOC-only subsystems.

Although some of the results could be derived by directly decoupling the associated ODEs of the system, the goal of applying CRNT extends beyond solving these particular systems. Our approach begins with relatively simple models as a foundational step, enabling a systematic application of CRNT techniques to more complex systems, particularly in scenarios where decoupling is infeasible or where parameter-free or parameter-minimal analyses are especially beneficial. CRNT offers a general framework to infer dynamic properties such as multistationarity and ACR without detailed knowledge of rate constants, making it a powerful tool for analyzing the qualitative behavior of reaction networks.

This framework can be expanded to include other carbon dioxide removal (CDR) strategies. Future work will explore necessary conditions for effective carbon reduction across various models, providing a way to evaluate and compare the long-term viability of different CDR approaches.

\appendix

\section{Details of showing the existence of positive steady states of the DOC system}
\label{details:existence}
The network $\NN$ can be decomposed into (stoichiometrically) independent subnetworks $\NN_1$ and $\NN_2$ given by 
\[\begin{aligned}
    \NN_1&:\quad
    A_1 + 2A_2 \xrightleftharpoons{\qquad} A_2 + 2A_1\\
    \NN_2&: 
    \quad \begin{tikzpicture}[node distance = {20mm}, baseline=(current  bounding  box.center)]
			\node (A17) {$A_{17}$};
            \node (A3) [above of=A17]{$A_3$};
            \node (A2) [left of=A3]{$A_2$};
            \node (A4) [left of=A17]{$A_4$};
            \draw[->] (A3) -- (A17);
            \draw[->] (A17) -- (A4);
            \draw[->] (A4) -- (A2);
            \draw [-left to] ($(A2.east) + (0pt, 2pt)$) -- ($(A3.west) + (0pt, 2pt)$);
    \draw [-left to] ($(A3.west) + (0pt, -2pt)$) -- ($(A2.east) + (0pt, -2pt)$);
        \end{tikzpicture}.
\end{aligned}.\] This decomposition is computed in Section \ref{subs:network:decomposition}, which can also be obtained using the MATLAB program in \cite{LubeniaINDECS}. In order to invoke the result of \cite{alamin2024positive}, we show that this decomposition also satisfies independence of its augmented matrix of kinetic order vectors (i.e. $\widehat T$-independence). Indeed, the $\widehat T$ matrix of the entire network $\NN$ is given by 

\[\widehat T = \begin{blockarray}{cccccccc}
        & \matindex{$A_1+2A_2$} & \matindex{$2A_1+A_2$} & \matindex{$A_2$} & \matindex{$A_3$} & \matindex{$A_4$} & \matindex{$A_{17}$} \\
        \begin{block}{c[ccccccc]}
        \matindex{$A_1$} & p_1 & p_2 & 0 & 0 & 0 & 0 \\
        \matindex{$A_2$} & q_1 & q_2 & 1 & 0 & 0 & 0 \\
        \matindex{$A_3$} & 0 & 0 & 0 & 1 & 0 & 0 \\
        \matindex{$A_4$} & 0 & 0 & 0 & 0 & 1 & 0 \\
        \matindex{$A_{17}$} & 0 & 0 & 0 & 0 & 0 & 1\\
        \matindex{$\NN_1$} & 1 & 1 & 0 & 0 & 0 & 0\\
        \matindex{$\NN_2$} & 0 & 0 & 1 & 1 & 1 & 1\\
        \end{block}
        \end{blockarray}\] which has rank six ($t = 6$) whenever $p_1 \neq p_2$ or $q_1\neq q_2.$ Moreover, if $p_1 = p_2$ and $q_1 = q_2,$ then the rank of $\widehat T$ is five ($t = 5$). The decomposition of $\NN$ into $\NN_1$ and $\NN_2$ gives rise to two $\widehat T$ matrices given by \[\widehat T_1 = \begin{blockarray}{ccc}
        & \matindex{$A_1+2A_2$} & \matindex{$2A_1+A_2$} \\
        \begin{block}{c[cc]}
        \matindex{$A_1$} & p_1 & p_2 \\
        \matindex{$A_2$} & q_1 & q_2 \\
        \matindex{$\NN_1$} & 1 & 1 \\
        \end{block}
        \end{blockarray}\quad\text{and}\quad \widehat T_2 = \begin{blockarray}{ccccc} & \matindex{$A_2$} & \matindex{$A_3$} & \matindex{$A_4$} & \matindex{$A_{17}$} \\
        \begin{block}{c[cccc]}
        \matindex{$A_2$} & 1 & 0 & 0 & 0 \\
        \matindex{$A_3$} & 0 & 1 & 0 & 0 \\
        \matindex{$A_4$} & 0 & 0 & 1 & 0 \\
        \matindex{$A_{17}$} & 0 & 0 & 0 & 1\\
        \matindex{$\NN_2$} & 1 & 1 & 1 & 1\\
        \end{block}
        \end{blockarray},\] which have ranks two ($t_1 = 2$) and four ($t_2 = 4$), respectively, whenever $p_1 \neq p_2$ or $q_1 \neq q_2.$ Now, if $p_1 = p_2$ and $q_1 = q_2$, then the rank of $\widehat T_1$ is one ($t_1 = 1$).  In any case, we get $\widehat T = \widehat T_1 \oplus \widehat T_2$ since their respective ranks add up, i.e. $t = t_1 + t_2.$ ($\widehat T$-independence). Invoking the result of Alamin and Hernandez \cite{alamin2024positive}, we conclude that the entire direct ocean capture system has positive steady states if and only if each subsystem induced by the independent subnetworks have positive steady states. This means that the existence of positive steady states of the DOC system, regardless of the values of $p_1, p_2, q_1,$ and $q_2$ may be determined through its independent subnetworks.

\section{Details of parametrization of positive steady states of the DOC system}
\label{details:parametrization}

To compute the positive steady state parametrization, we follow the steps provided in \cite{hernandez2023network,hernandez2023framework}. The first step is to get the finest independent decomposition of the whole network $\NN$ as computed in Section \ref{subs:network:decomposition}. Next, we get the positive steady states of each subnetwork ($\NN_1$ and $\NN_2$) individually.

\subsection{Computation of positive steady states of \texorpdfstring{$\NN_1$}{\NN1}}

The following steps are due to Johnston et al. \cite{johnston2019deficiency} via the so-called ``network translation.'' If we can find such a network that is weakly reversible and deficiency zero, then we can compute the positive steady states via this method. For a more detailed discussion of the method, please refer to \cite{hernandez2023network,johnston2019deficiency}.

\noindent STEP 1: Find a weakly reversible and deficiency zero translated network. Translating a network can be done by adding or subtracting the same term to both sides of the reactions (preserving the stoichiometric matrix of the network) but considering the original kinetic vectors (preserving the rate functions).

\begin{table}[h!]
    \centering
    \begin{tabular}{c c c}
    \textbf{Original Network} & \textbf{Kinetic Order Vector} \\
    $R_1: A_1 + 2A_2 \longrightarrow 2A_1 + A_2$ & $p_1A_1 + q_1A_2=[p_1,q_1]^{\top}$ \\
    $R_2: 2A_1 + A_2 \longrightarrow A_1 + 2A_2$ & $p_2A_1 + q_2A_2=[p_2,q_2]^{\top}$ 
    \end{tabular}
\end{table} 

\begin{center}
\textbf{Translated Network}
\end{center}

\begin{center}
\begin{tikzpicture}
        \tikzset{vertex/.style ={rectangle, draw, minimum width =10pt,{minimum size=2em}}}
        \tikzset{edge/.style = {->,> = {Stealth[length=2mm, width=2mm]}}}
        \node[vertex,text width=3cm,text centered] (A) at (0,0) {{ \circled{1}  \\ $A_2$ \\ $(p_1A_1 + q_1A_2)$}};
        \node[vertex,text width=3cm,text centered] (AE) at (6,0) {{ \circled{2} \\ $A_1$ \\ $(p_2A_1 + q_2A_2)$}};
        \draw[edge, thick] (A.11) to (AE.169)
        node[above,xshift=-12mm] {$k_1$};
        \draw[edge, thick] (AE.190) to (A.350)
        node[below,xshift=15mm] {$k_2$};
\end{tikzpicture}
\end{center}


\noindent STEP 2: Get all the spanning trees, with edges labeled by rate constants, towards each node.

\begin{table}[h!]
    \centering
    \begin{tabular}{c c}
    \textbf{Towards 1} & \textbf{Towards 2} \\
    $k_2: 2 \longrightarrow 1$ & $k_1: 1 \longrightarrow 2$ \\
    $K_1 = k_2$ & $K_2 = k_1$ 
    \end{tabular}
\end{table} 

\noindent STEP 3: Choose any spanning tree containing all the nodes. (Here, we choose $1 \to 2$.)
Furthermore, we compute $\kappa_{i\to i'}=\dfrac{K_{i'}}{K_i}$ and get the kinetic difference(s) (i.e., difference between the kinetic vectors given inside the parentheses in the translated network) associated to the edge(s) of the tree.

\begin{table}[h!]
    \centering
    \begin{tabular}{c|c}
     $\kappa_{1 \rightarrow 2} = \dfrac{K_2}{K_1} = \dfrac{k_1}{k_2}$ & $(p_2 - p_1)A_1 + (q_2 - q_1)A_2$ 
    \end{tabular}
\end{table} 

\noindent STEP 4: Compute matrices $M$, $H$, and $B$.
We have $M = \begin{bmatrix} p_2 - p_1 & q_2 - q_1 \end{bmatrix}$ (the matrix of kinetic differences). We find $H = \begin{bmatrix} h_1 & h_2\end{bmatrix}^\top$ such that $MHM = M$, i.e.,
$$
\begin{bmatrix} p_2 - p_1 & q_2 - q_1 \end{bmatrix}
\begin{bmatrix} h_1 \\ h_2
\end{bmatrix}
\begin{bmatrix} p_2 - p_1 & q_2 - q_1 \end{bmatrix} = \begin{bmatrix} p_2 - p_1 & q_2 - q_1 \end{bmatrix} 
$$
We have $H = \begin{bmatrix} \dfrac{1}{p_2 - p_1} \\ 0 \end{bmatrix}$. \\

Let $B = \begin{bmatrix} b_1 & b_2\end{bmatrix}^\top$. We find matrix $B$ such that $\text{ker} \hspace{0.1cm} M = B$. We have

\begin{align*}
\begin{bmatrix} p_2 - p_1 & q_2 - q_1 \end{bmatrix} \begin{bmatrix} b_1 \\ b_2 \end{bmatrix} &= 0 \\
(p_2 - p_1)b_1 + (q_2 - q_1)b_2 &= 0 \\
(q_2 - q_1)b_2 &= (p_1 - p_2)b_1 \\
b_2 &= \dfrac{(p_1 - p_2)b_1}{q_2 - q_1}.
\end{align*}
If $b_1 = q_2 - q_1$, then $b_2 = p_1 - p_2$. So $B = \begin{bmatrix} q_2 - q_1 \\ p_1 - p_2\end{bmatrix}$. \\

\noindent STEP 5: Establish positive steady states. The values of $a_1$ and $a_2$ using the entries of matrices $H$ and $B$ as exponents are
\begin{align*}
   a_1 &= (\kappa_{2 \rightarrow 1})^\frac{1}{p_2 - p_1} \hspace{0.1cm} \tau^{q_2 - q_1} = \left(\dfrac{k_1}{k_2}\right)^\frac{1}{p_2 - p_1} \hspace{0.1cm} \tau^{q_2 - q_1}\\
   a_2 &= (\kappa_{2 \rightarrow 1})^0 \hspace{0.1cm} \tau^{p_1 - p_2} = \tau^{p_1 - p_2}
\end{align*}
with $\tau>0$, a free parameter (only one) since the matrix $B$ only has one column vector. This covers the case when $p_1 \ne p_2$.

However, we can also choose $H$ to be $\begin{bmatrix} 0 \\ \dfrac{1}{q_2 - q_1} \end{bmatrix}$. In this case, the parametrization is
\begin{align*}
   a_1 &= (\kappa_{2 \rightarrow 1})^0 \hspace{0.1cm} \tau^{q_2 - q_1} = \tau^{q_2 - q_1}\\
   a_2 &= (\kappa_{2 \rightarrow 1})^\frac{1}{q_2 - q_1} \hspace{0.1cm} \tau^{p_1 - p_2} = \left(\dfrac{k_1}{k_2}\right)^\frac{1}{q_2 - q_1} \hspace{0.1cm} \tau^{p_1 - p_2}.
\end{align*}
This covers the case when $q_1 \ne q_2$.

\subsection{Computation of positive steady states of \texorpdfstring{$\NN_2$}{\NN2}}

Here, we provide a parametrization of $\mathcal{N}_2$ via the computational package COMPILES (COMPutIng anaLytic stEady States) developed in \cite{hernandez2023framework}, which is built in MATLAB. It derives a steady state parametrization of the network by decomposing the CRN into independent subnetworks and combines parametrizations of the subnetworks. Note that COMPILES is only applicable for mass action systems.\\

\noindent \textbf{Code} (used on the script file)

\begin{verbatim}
model.id = 'Direct Ocean Capture';
model = addReaction(model, 'A1+2A2<->2A1+A2', ...               
                           {'A1', 'A2'}, {1, 2}, [1, 2], ...    
                           {'A1', 'A2'}, {2, 1}, [2, 1], ...    
                           true);                               
model = addReaction(model, 'A2<->A3', ...
                           {'A2'}, {1}, [1], ...
                           {'A3'}, {1}, [1], ...
                           true);
model = addReaction(model, 'A4->A2', ...
                           {'A4'}, {1}, [1], ...
                           {'A2'}, {1}, [ ], ...
                           false);
model = addReaction(model, 'A17->A4', ...
                           {'A17'}, {1}, [1], ...
                           {'A4'}, {1}, [ ], ...
                           false);
model = addReaction(model, 'A3->A17', ...
                           {'A3'}, {1}, [1], ...
                           {'A17'}, {1}, [ ], ...
                           false);

[equation, species, free_parameter, conservation_law, model] 
= steadyState(model); 
\end{verbatim} 

\noindent \textbf{Output}

\begin{verbatim}
The network has 2 subnetworks.

- Subnetwork 1 -

R1: A1+2A2->2A1+A2
R2: 2A1+A2->A1+2A2

Solving Subnetwork 1...

A1 = (k1*tau1)/k2 
A2 = tau1 

- Subnetwork 2 -

R3: A2->A3
R4: A3->A2
R5: A4->A2
R6: A17->A4
R7: A3->A17

Solving Subnetwork 2...

A17 = (k5*tau2)/k6 
A2 = (k5*tau2*(k4 + k7))/(k3*k7) 
A3 = (k5*tau2)/k7 
A4 = tau2 

Solving positive steady state parametrization of the entire 
network...

The solution is as follows.

A1 = (A4*k1*k5*(k4 + k7))/(k2*k3*k7)
A2 = (A4*k5*(k4 + k7))/(k3*k7)
A3 = (A4*k5)/k7
A17 = (A4*k5)/k6
Free parameter: A4
\end{verbatim}

We focus solely on the solution for the second subnetwork in the output, as it follows the mass action formalism, whereas the first network follows the power law formalism and was computed earlier in the previous subsection.

Hence the obtained parametrized steady state solution for $\NN_2$ is given by
\begin{align*}
    \begin{cases}
    a_2 = \dfrac{(k_5 \omega)(k_4 + k_7)}{k_3 k_7} \\
    a_3 = \dfrac{k_5 \omega}{k_7} \\
    a_4 = \omega \\
    a_{17} = \dfrac{k_5 \omega}{k_6}
    \end{cases}.
\end{align*}

\subsection{Computation of positive steady states of the DOC system}

First, we consider the case when $p_1 \ne p_2$, we merge the obtained positive steady states in the preceding two subsections (the values of $a_2$ which is common to both subnetworks must agree) to obtain the following steady state parametrization of the whole network:
\begin{align*}
    \begin{cases}
    a_1 = \left(\dfrac{k_1}{k_2}\right)^{\frac{1}{p_2-p_1}} \tau^{q_2-q_1} \\
    a_2 = \tau^{p_1-p_2} \\
    a_3 = \dfrac{k_3}{k_4 + k_7} \tau^{p_1-p_2} \\
    a_4 = \dfrac{k_3k_7}{k_5(k_4+k_7)} \tau^{p_1-p_2} \\
    a_{17} = \dfrac{k_3k_7}{k_6(k_4+k_7)} \tau^{p_1-p_2} \\
    \text{Free parameter:} \hspace{0.2cm} \tau > 0.
    \end{cases}
\end{align*}

Recall that the ODEs for $\mathcal{N}$ are
    \begin{align*}
    \dfrac{da_1}{dt} &= k_1a_1^{p_1}a_2^{q_1} - k_2a_1^{p_2}a_2^{q_2}  \\
    \dfrac{da_2}{dt} &= k_2a_1^{p_2}a_2^{q_2} - k_1a_1^{p_1}a_2^{q_1} - k_3a_2 + k_4a_3 + k_5a_4  \\
    \dfrac{da_3}{dt} &= k_3a_2 - k_4a_3 - k_7a_3  \\
    \dfrac{da_4}{dt} &= k_6a_{17} - k_5a_4 \\
    \dfrac{da_{17}}{dt} &= k_7a_3 - k_6a_{17}.
\end{align*}
We substitute the obtained parameterized steady-state solution into each of the ODEs and verify that it indeed makes the right-hand side of each equation in the ODE system for the entire network $\mathcal{N}$ equal to zero.

Second, we consider the case when $q_1 \ne q_2$. Following the same steps from the previous case, we arrive at the parametrization
\begin{align*}
    \begin{cases}
    a_1 = \tau^{q_2-q_1} \\
    a_2 = \left(\dfrac{k_1}{k_2}\right)^{\frac{1}{q_2-q_1}} \tau^{p_1-p_2} \\
    a_3 = \dfrac{k_3}{k_4 + k_7} \left(\dfrac{k_1}{k_2}\right)^{\frac{1}{q_2-q_1}} \tau^{p_1-p_2} \\
    a_4 = \dfrac{k_3k_7}{k_5(k_4+k_7)} \left(\dfrac{k_1}{k_2}\right)^{\frac{1}{q_2-q_1}} \tau^{p_1-p_2} \\
    a_{17} = \dfrac{k_3k_7}{k_6(k_4+k_7)} \left(\dfrac{k_1}{k_2}\right)^{\frac{1}{q_2-q_1}} \tau^{p_1-p_2} \\
    \text{Free parameter:} \hspace{0.2cm} \tau > 0.
    \end{cases}
\end{align*}

\section{Details of the analysis for the conditions of multistationarity in DOC systems}
\label{details:multistationarity}

To determine some sufficient conditions for the direct ocean capture model to admit multiple steady states, we utilize the following result by M\"uller and Regensburger \cite{regensburger}.

\begin{theorem} \label{muller}
    If for a weakly reversible generalized mass action system with ${\rm sign} (S) \cap {\rm sign} (\tilde{S})^{\perp} \neq \{0\}$, then there is a stoichiometric class with more than one (complex balanced) steady state.
\end{theorem}

The theorem tells us that for weakly reversible generalized mass action systems, a sufficient condition for the system to be multistationary is the existence of a non-trivial vector whose sign pattern is the same as that of the stoichiometric subspace $S$ and the orthogonal complement of kinetic flux subspace $\tilde S.$

First, we solve for the sign pattern of $\tilde S.$ Note that $\tilde{S} = \text{Im}\,(\tilde{Y} \cdot I_a)$ where
\begin{center}
    $\tilde{Y}=\begin{blockarray}{cccccccc}
        & \matindex{$A_1+2A_2$} & \matindex{$2A_1+A_2$} & \matindex{$A_2$} & \matindex{$A_3$} & \matindex{$A_4$} & \matindex{$A_{17}$} \\
        \begin{block}{c[ccccccc]}
        \matindex{$A_1$} & p_1 & p_2 & 0 & 0 & 0 & 0 \\
        \matindex{$A_2$} & q_1 & q_2 & 1 & 0 & 0 & 0 \\
        \matindex{$A_3$} & 0 & 0 & 0 & 1 & 0 & 0 \\
        \matindex{$A_4$} & 0 & 0 & 0 & 0 & 1 & 0 \\
        \matindex{$A_{17}$} & 0 & 0 & 0 & 0 & 0 & 1 \\
        \end{block}
        \end{blockarray}$
\end{center}
and
\begin{center}
    $I_a=\begin{blockarray}{cccccccc}
        & \matindex{$R_1$} & \matindex{$R_2$} & \matindex{$R_3$} & \matindex{$R_4$} & \matindex{$R_5$} & \matindex{$R_6$} & \matindex{$R_7$} \\
        \begin{block}{c[ccccccc]}
        \matindex{$A_1+2A_2$} & -1 & 1 & 0 & 0 & 0 & 0 & 0 \\
        \matindex{$2A_1+A_2$} & 1 & -1 & 0 & 0 & 0 & 0 & 0 \\
        \matindex{$A_2$} & 0 & 0 & -1 & 1 & 1 & 0 & 0 \\
        \matindex{$A_3$} & 0 & 0 & 1 & -1 & 0 & 0 & -1 \\
        \matindex{$A_4$} & 0 & 0 & 0 & 0 & -1 & 1 & 0 \\
        \matindex{$A_{17}$} & 0 & 0 & 0 & 0 & 0 & -1 & 1 \\
        \end{block}
        \end{blockarray}$.
\end{center} Here, the $\tilde Y$ matrix is defined using the kinetic order vectors of the system (see \cite{regensburger}) and $I_a$ is the incidence matrix of the network. Hence, 
$$\tilde{Y} \cdot I_a = \begin{bmatrix} 
p_2-p_1 & p_1-p_2 & 0 & 0 & 0 & 0 & 0 \\
q_2-q_1 & q_1-q_2 & -1 & 1 & 1 & 0 & 0 \\
0 & 0 & 1 & -1 & 0 & 0 & -1 \\
0 & 0 & 0 & 0 & -1 & 1 & 0 \\
0 & 0 & 0 & 0 & 0 & -1 & 1
\end{bmatrix}$$
$$\Rightarrow \tilde{S} = \text{Im} \hspace{0.1cm} (\tilde{Y} \cdot I_a) = \sf{span} \hspace{0.1cm} \left \{\begin{bmatrix} p_2-p_1 \\ q_2-q_1 \\ 0 \\ 0 \\ 0 \end{bmatrix}, \begin{bmatrix} 0 \\ -1 \\ 1 \\ 0 \\ 0 \end{bmatrix}, \begin{bmatrix} 0 \\ 0 \\ -1 \\ 0 \\ 1 \end{bmatrix}, \begin{bmatrix} 0 \\ 0 \\ 0 \\ 1 \\ -1 \end{bmatrix} \right \}.$$
The orthogonal complement $(\tilde{S})^{\perp}$ of $\tilde{S}$ is given by
$$(\tilde{S})^{\perp} = \rm{span} \left \{\begin{bmatrix} \frac{q_1-q_2}{p_2-p_1} \\ 1 \\ 1 \\ 1 \\ 1 \end{bmatrix} \right \} = \rm{span} 
 \left \{\begin{bmatrix} -Q \\ 1 \\ 1 \\ 1 \\ 1 \end{bmatrix} \right \} = \rm{span}  \left \{\begin{bmatrix} -1 \\ {R} \\ {R} \\ {R} \\ {R} \end{bmatrix} \right \}$$
where ${R} = \dfrac{p_2-p_1}{q_2-q_1}$ and ${Q} = \dfrac{q_2-q_1}{p_2-p_1},$ as defined.

We now investigate the multiplicity of steady states for positive ($R > 0$), negative ($R < 0)$, $P$-null ($R = 0$ and defined), and $Q$-null ($Q = 0$ and defined) systems. 

First, for positive DOC systems, i.e., $R > 0\, (Q > 0),$ we have \[{\rm sign}(\tilde S^\perp) =  \left \{\begin{bmatrix} - \\ + \\ + \\ + \\ + \end{bmatrix}, \begin{bmatrix} + \\ - \\ - \\ - \\ - \end{bmatrix} \right \}.\] 

Indeed, if we let $x$ be in the stoichiometric subspace $S$ given by \[S = \text{span }\left \{\begin{bmatrix} 1 \\ -1 \\ 0 \\ 0 \\ 0 \end{bmatrix}, \begin{bmatrix} 0 \\ -1 \\ 1 \\ 0 \\ 0 \end{bmatrix}, \begin{bmatrix} 0 \\ 0 \\ -1 \\ 0 \\ 1 \end{bmatrix}, \begin{bmatrix} 0 \\ 0 \\ 0 \\ 1 \\ -1 \end{bmatrix} \right \},\] then \[x = a_1 \begin{bmatrix} 1 \\ -1 \\ 0 \\ 0 \\ 0 \end{bmatrix} + a_2 \begin{bmatrix} 0 \\ -1 \\ 1 \\ 0 \\ 0 \end{bmatrix} + a_3 \begin{bmatrix} 0 \\ 0 \\ -1 \\ 0 \\ 1 \end{bmatrix} + a_4 \begin{bmatrix} 0 \\ 0 \\ 0 \\ 1 \\ -1 \end{bmatrix} = \begin{bmatrix} a_1 \\ -a_1-a_2 \\ a_2-a_3 \\ a_4 \\ a_3-a_4 \end{bmatrix}.\] We can then choose $a_1 > 0$ and $a_2 < a_3 < a_4 < 0$ so that we have \[{\rm sign}(x) = \begin{bmatrix} + \\ - \\ - \\ - \\ - \end{bmatrix} \in {\rm sign}(\tilde S^\perp)\] and thus ${\rm sign}(x) \cap {\rm sign}(\tilde{S})^{\perp} \neq \{0\}.$ Therefore, by Theorem \ref{muller}, any positive DOC system is multistationary.

Now, for negative DOC systems, we cannot utilize Theorem \ref{muller} to conclude monostationarity. Because of this, we employ a different criterion to conclude when the system is monostationary. The following computational method introduced by Wiuf and Feliu \cite{wiuf, interacting} reveals network injectivity for a specific subset of the collection of negative DOC systems.

\begin{theorem} (\textit{Feliu and Wiuf, 2013} \cite{interacting}) \label{powerlawinjectivity}
    The interaction network with power law kinetics and fixed kinetic orders is injective if and only if the determinant of $M^*$ is a nonzero homogeneous polynomial with all coefficients being positive or all being negative.
\end{theorem}

Since network injectivity implies monostationarity \cite{feinberg2019crnt}, we can study the individual signs of $p_1, p_2, q_1,$ and $q_2$ to know when the system achieve monostationarity.

The matrix $M^*$ in the theorem is defined using the kinetic order matrix $F$ and stoichiometric matrix $N$ of the network. Indeed, we consider symbolic vectors $k = (k_1, \ldots, k_m)$ and $z = (z_1, \ldots, z_r)$ and define $M = N \text{diag} (z) F \text{diag} (k)$. Taking $\{w_1, \ldots, w_d\}$ to be a basis of the left kernel of $N$ and $i_1, \ldots, i_d$ row indices as above, we can write $i_j$ to denote the index of the first nonzero entry of $w^j$ \cite{rau}. From this, we define the $m \times m$ matrix $M^*$, by replacing the $i_j$-th row of $M$ by $w_j$. Note that matrix $M^*$ is a symbolic matrix in $z_*$ and $k_*$.

The stoichiometric matrix and matrix of kinetic order vectors for our DOC system is given by 
\begin{align*}
N &= \begin{bmatrix}
    1 & -1 & 0 & 0 & 0 & 0 & 0 \\
    -1 & 1 & -1 & 1 & 1 & 0 & 0 \\
    0 & 0 & 1 & -1 & 0 & 0 & -1 \\
    0 & 0 & 0 & 0 & -1 & 1 & 0 \\
    0 & 0 & 0 & 0 & 0 & -1 & 1 
\end{bmatrix}, \quad \text{and} \\
F &= \begin{bmatrix}
    p_1 & q_1 & 0 & 0 & 0 \\
    p_2 & q_2 & 0 & 0 & 0 \\
    0 & 1 & 0 & 0 & 0 \\
    0 & 0 & 1 & 0 & 0 \\
    0 & 0 & 0 & 1 & 0 \\
    0 & 0 & 0 & 0 & 1 \\
    0 & 0 & 1 & 0 & 0
\end{bmatrix}.
\end{align*} respectively. Moreover, given the symbolic vectors $k$ and $z,$ we have 
    \begin{align*}
    \text{diag}(z) &= \begin{bmatrix}
        z_1 & 0 & 0 & 0 & 0 & 0 & 0 \\
        0 & z_2 & 0 & 0 & 0 & 0 & 0 \\
        0 & 0 & z_3 & 0 & 0 & 0 & 0 \\
        0 & 0 & 0 & z_4 & 0 & 0 & 0 \\
        0 & 0 & 0 & 0 & z_5 & 0 & 0 \\
        0 & 0 & 0 & 0 & 0 & z_6 & 0 \\
        0 & 0 & 0 & 0 & 0 & 0 & z_7
    \end{bmatrix}, \quad \text{and} \\
    \text{diag}(k) &= \begin{bmatrix}
        k_1 & 0 & 0 & 0 & 0 \\
        0 & k_2 & 0 & 0 & 0 \\
        0 & 0 & k_3 & 0 & 0 \\
        0 & 0 & 0 & k_4 & 0 \\
        0 & 0 & 0 & 0 & k_5
    \end{bmatrix}.
    \end{align*}

We construct the matrix $M = N\text{diag}(z)F\text{diag}(k)$ and obtain
\begin{center}
\resizebox{\textwidth}{!}{$
M = \begin{bmatrix}
    k_1p_1z_1 - k_1p_2z_2 & k_2q_1z_1 - k_2q_2z_2 & 0 & 0 & 0 \\
    -k_1p_1z_1 + k_1p_2z_2 & -k_2q_1z_1 + k_2q_2z_2 - k_2z_3 & k_3z_4 & k_4z_5 & 0 \\
    0 & k_2z_3 & -k_3z_4 - k_3z_7 & 0 & 0 \\
    0 & 0 & 0 & -k_4z_5 & k_5z_6 \\
    0 & 0 & k_3z_7 & 0 & -k_5z_6
\end{bmatrix}
$}
\end{center}
Now, the basis of the left kernel of $N$ is $\{[1, 1, 1, 1, 1]\}$. This row vector will replace the first row of the matrix $M$. Therefore, we have our matrix $M^*$ given by
\begin{center}
\resizebox{\textwidth}{!}{$
M^* = \begin{bmatrix}
    1 & 1 & 1 & 1 & 1 \\
    -k_1p_1z_1 + k_1p_2z_2 & -k_2q_1z_1 + k_2q_2z_2 - k_2z_3 & k_3z_4 & k_4z_5 & 0 \\
    0 & k_2z_3 & -k_3z_4 - k_3z_7 & 0 & 0 \\
    0 & 0 & 0 & -k_4z_5 & k_5z_6 \\
    0 & 0 & k_3z_7 & 0 & -k_5z_6
\end{bmatrix}
$}
\end{center}
The determinant of $M^*$, using computer software MATLAB is found to be 
\begin{align*}
    \det M^* = &- \underline{p_1}k_1k_2k_4k_5z_1z_3z_5z_6 - \underline{p_1}k_1k_2k_3k_4z_1z_3z_5z_7 \\ &-	
    \underline{p_1}k_1k_2k_3k_5z_1z_3z_6z_7 - \underline{p_1}k_1k_3k_4k_5z_1z_4z_5z_6 \\ &-	
    \underline{p_1}k_1k_3k_4k_5z_1z_5z_6z_7 + \underline{p_2}k_1k_2k_3k_4z_2z_3z_5z_7 \\ &+ 	
    \underline{p_2}k_1k_2k_4k_5z_2z_3z_5z_6 + \underline{p_2}k_1k_2k_3k_5z_2z_3z_6z_7 \\&+\underline{p_2}k_1k_3k_4k_5z_2z_4z_5z_6 + \underline{p_2}k_1k_3k_4k_5z_2z_5z_6z_7 \\ &+ \underline{q_1}k_2k_3k_4k_5z_1z_4z_5z_6 + \underline{q_1}k_2k_3k_4k_5z_1z_5z_6z_7 \\ &-\underline{q_2}k_2k_3k_4k_5z_2z_4z_5z_6 - \underline{q_2}k_2k_3k_4k_5z_2z_5z_6z_7.
\end{align*}

Hence, for $p_1 < 0, p_2 > 0, q_1 > 0$, and $q_2 < 0$, all the terms of the determinant are positive, and for $p_1 > 0, p_2 < 0, q_1 < 0$, and $q_2 > 0$, all the terms of the determinant are negative. By Theorem \ref{powerlawinjectivity}, the systems in these cases are injective, and hence monostationary. These conditions, although sufficient, are not necessary for monostationarity. Therefore, only a subset of our negative DOC systems, specifically systems satisfying either (i) $p_1, q_2 > 0$ and $p_2, q_1 < 0$ or (ii) $p_1, q_2 < 0$ and $p_2, q_1 > 0$, exhibit monostationarity.

Finally, for the $P$-null and $Q$-null DOC systems, we investigate their induced ODEs and arrive at the following theorems:

\begin{theorem}
\label{thm-monoP0}
    All ${\sf DOC}_{P_0}$ systems are monostationary.
    \begin{proof}
    Suppose otherwise and let $E_1$ and $E_2$ be two distinct equilibria in the same stoichiometric class of a $P$-null system. Since we achieve ACR on $A_2, A_3, A_4,$ and $A_{17},$ the concentrations at steady state for such species are fixed for any set of rate constants. Hence, $E_1$ and $E_2$ differ in their concentration of $A_1.$ Note that by the induced ODEs of DOC systems, we have the conservation law \[A_1'(t) + A_2'(t) + A_3'(t) + A_4'(t) + A_{17}'(t) = 0.\] Thus, at its steady state, we have \[a_1 = A_0 - a_2 - a_3 - a_4 - a_{17}\] where $A_0$ is the initial concentration of carbon in the system, which remains fixed in the same set of rate constants. From this, we have \[a_1 = A_0 - \left(\frac{k_1}{k_2}\right)^{\frac{1}{q_2 - q_1}}\left(1  + \frac{k_3}{k_4 + {k_7}} + \frac{k_3{k_7}}{{k_5}(k_4 + {k_7})} + \frac{k_3{k_7}}{k_6(k_4 + {k_7})}\right).\] Thus, noting that $A_0$ remains fixed, the concentration of $a_1$ is unique at steady state for any fixed set of rate constants. Therefore, the system is monostationary.
    \end{proof}
\end{theorem}

A similar approach may be done to conclude the monostationarity of $Q$-null DOC systems. 

\begin{theorem}
\label{thm-monoQ0}
    All ${\sf DOC}_{Q_0}$ systems are monostationary.
    \begin{proof}
        Following the proof of Theorem \ref{thm-monoP0}, we utilize the same conservation law so that at steady state, we also have \[a_1 = A_0 - a_2 - a_3 - a_4 - a_{17}\] where $A_0$ is the initial concentration of carbon in the system. Using the parametrization of the steady states of systems in ${\sf DOC}_{Q_0},$ we have \[a_1 = A_0 - \tau^{p_1 - p_2}\left(1 + \frac{k_3}{k_4 + k_7} + \frac{k_3k_7}{k_5(k_4 + k_7)} + \frac{k_3k_7}{k_6(k_4 + k_7)} \right)\] where $\tau > 0.$ Since we achieve ACR on $A_1$ only for $Q$-null DOC systems, the free parameter $\tau$ becomes determined since $a_1$ is the same across all sets of rate constants, i.e. \[\tau = \left(\frac{A_0 - a_1}{1 + \frac{k_3}{k_4 + k_7} + \frac{k_3k_7}{k_5(k_4 + k_7)} + \frac{k_3k_7}{k_6(k_4 + k_7)}}\right)^\frac{1}{p_1 - p_2}.\] This implies then that the concentrations of $A_2, A_3, A_4,$ and $A_{17}$ at steady state is unique. Therefore, systems in ${\sf DOC}_{Q_0}$ are monostationary.
    \end{proof}
\end{theorem}

\section{DOC admits a positive steady state for every stoichiometric class}
\label{details:stoichiometric:class}

We can observe in the ODE system of the DOC in Section 2.1 that the following equation holds $$\dfrac{da_1}{dt}+\dfrac{da_2}{dt}+\dfrac{da_3}{dt}+\dfrac{da_4}{dt}+\dfrac{da_{17}}{dt}=0$$
by adding the right hand side of all the equations in the ODE system.
This means that the total concentration is constant for any time $t$, i.e., $a_1+a_2+a_3+a_4+a_{17}=T>0$, the conservation equation. In particular, at the positive steady steady state, 
\begin{center}
\resizebox{\textwidth}{!}{%
$\begin{aligned}
\left(\dfrac{k_1}{k_2}\right)^{\frac{1}{p_2-p_1}} \tau^{q_2-q_1} + \tau^{p_1-p_2} + \dfrac{k_3}{k_4 + k_7} \tau^{p_1-p_2} + \dfrac{k_3k_7}{k_5(k_4+k_7)} \tau^{p_1-p_2} + \dfrac{k_3k_7}{k_6(k_4+k_7)} \tau^{p_1-p_2} &= T \\
\left(\dfrac{k_1}{k_2}\right)^{\frac{1}{p_2-p_1}} \tau^{q_2-q_1} + \left[ 1 + \dfrac{k_3}{k_4 + k_7}+\dfrac{k_3k_7}{k_5(k_4+k_7)}+\dfrac{k_3k_7}{k_6(k_4+k_7)}\right]\tau^{p_1-p_2} &= T
\end{aligned}$
}
\end{center}
by replacing the concentrations using the parametrization of positive steady states computed for the DOC system.
Thus, the equation has the form $a\tau^{q_2-q_1}+b\tau^{p_1-p_2}-c=0.$ where $a,b,c>0$.

We now analyze the existence of solutions to the equation by examining sign changes, based on a generalization of Descartes' rule of signs for counting the number of positive solutions \cite{zeros}. There are nine possible combinations of the values of the exponents $q_2-q_1$ and $p_1-p_2$, since each exponent difference can be positive, negative, or zero. The case where both differences are zero is not included. Furthermore, other combinations can be combined into a single condition. Arranging the generalized polynomial so that the exponents are listed in decreasing order gives rise to the following scenarios:
\begin{enumerate}
    \item $q_2-q_1>0$ and $p_1-p_2>0$: one sign change
    \item $q_2-q_1<0$ and $p_1-p_2<0$: one sign change
    \item $q_2-q_1=0$ and $p_1-p_2\ne 0$: one sign change
    \item $q_2-q_1\ne 0$ and $p_1-p_2=0$: one sign change
    \item $q_2-q_1> 0$ and $p_1-p_2<0$: two sign changes
    \item $q_2-q_1< 0$ and $p_1-p_2>0$: two sign changes
\end{enumerate}
For the first four cases, there is exactly one positive solution. In the last two cases, the number of positive solutions is either zero or two. However, we have demonstrated that the system exhibits multistationarity in these two cases. Thus, a positive steady state exists for each stoichiometric class in any DOC system.

\section{Simulations confirming monostationarity or multistationarity of DOC systems}
\label{appendix:validation:multi}

In this section, we validate our results on multistationarity for the DOC systems summarized in Table~\ref{tab:summary} of the main manuscript by plotting the functions
\[
y(\tau) = C_1 \tau^{q_2 - q_1} + C_2 \tau^{p_1 - p_2} - T
\]
for non-P-null systems (i.e., positive, negative, and Q-null DOC systems), and
\[
z(\tau) = \tau^{q_2 - q_1} + C_2 C_3 \tau^{p_1 - p_2} - T
\]
for P-null systems, as derived in Appendix~\ref{details:stoichiometric:class}. Here, \( T \) represents the total concentration given by the conservation relation.

We consider four sets of kinetic orders \( (p_1, q_1, p_2, q_2) \), each corresponding to a different DOC type:
\begin{enumerate}
    \item \( (1.5, 1.0, 2.5, 3.0) \) for positive DOC  (plotted using the function \( y(\tau) \));
    \item \( (-1.0, 1.5, 1.0, -1.5) \) for negative DOC (plotted using the function \( y(\tau) \));
    \item \( (1.0, 0.5, 1.0, 2.5) \) for P-null DOC (plotted using the function \( z(\tau) \));
    \item \( (3.0, 1.5, 1.0, 1.5) \) for Q-null DOC (plotted using the function \( y(\tau) \)).
\end{enumerate}

The asterisks on the graphs indicate the values of \( \tau \) at which the corresponding function crosses the x-axis. These points identify the steady state solutions of the system. If the function intersects the x-axis at two points, the system admits two steady states (multistationarity) while a single intersection implies monostationarity. Each root of \( \tau \) corresponds to a distinct set of steady-state concentrations via the parametrization formulas.

\begin{figure}[H]
    \begin{center}
    \includegraphics[width=11cm,height=7cm,keepaspectratio]{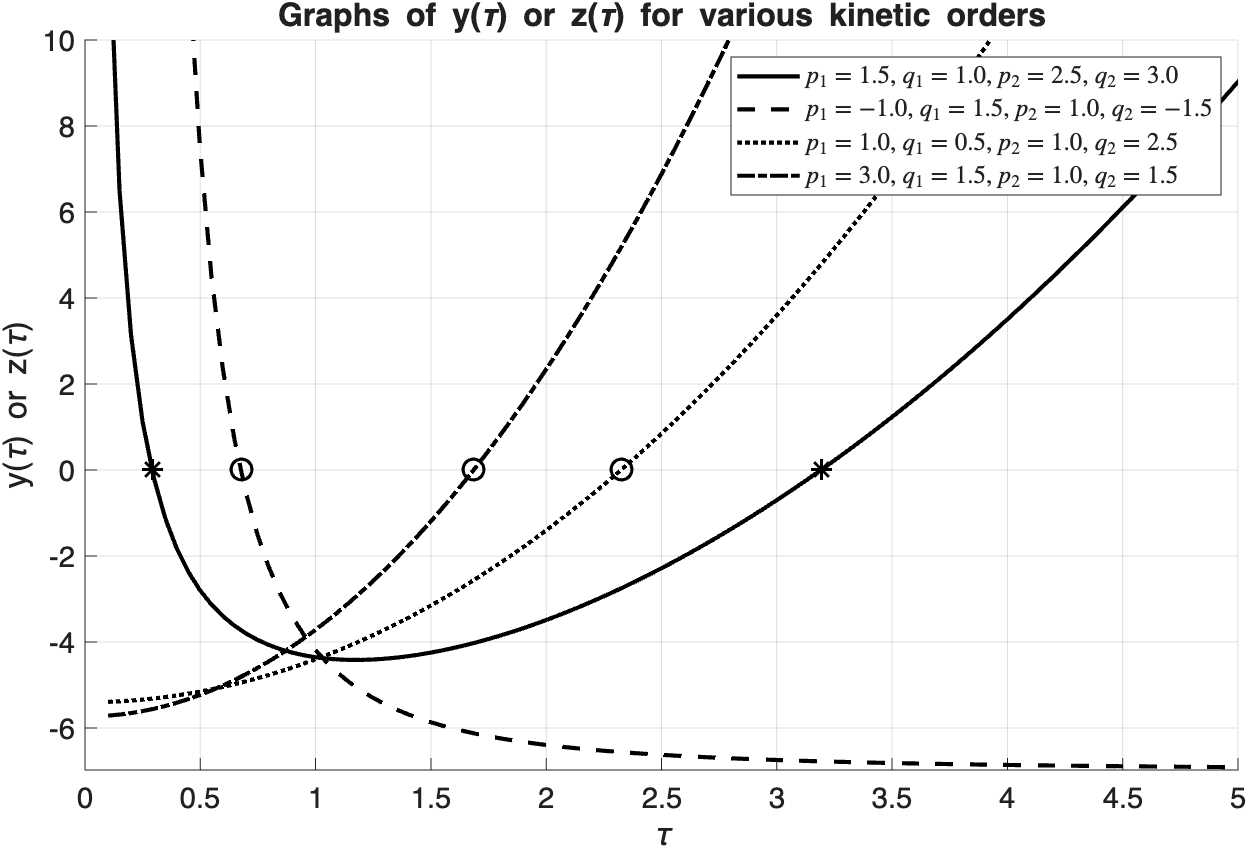}
        \caption{Plots of the the functions $y(\tau) = C_1 \tau^{q_2 - q_1} + C_2 \tau^{p_1 - p_2} - T$
    for non-P-null systems (positive, negative, and Q-null DOC systems),
    and $z(\tau) = \tau^{q_2 - q_1} + C_2 C_3 \tau^{p_1 - p_2} - T$
    for P-null systems,
    where $T$ is the total concentration from the conservation equation. These expressions result from substituting the steady state parametrizations into the conservation equation. Figure \ref{fig:MultistationaritySimulation} shows plots of the functions with four sets of kinetic orders $(p_1, q_1, p_2, q_2)$ correspond to various DOC types: 1. $(1.5, 1.0, 2.5, 3.0)$ corresponds to a positive DOC, 2. $(-1.0, 1.5, 1.0, -1.5)$ corresponds to a negative DOC, 3. $(1.0, 0.5, 1.0, 2.5)$ is for a P-null DOC and uses the function $z(\tau)$, and 4. $(3.0, 1.5, 1.0, 1.5)$ corresponds to a Q-null DOC.
Asterisks on the plot mark the $\tau$ values where the functions cross the x-axis, indicating monostationarity (with one solution) or multistationarity (with two solutions) of the DOC systems.
}
    \label{fig:MultistationaritySimulation}
    \end{center}
    \end{figure}

\section{Simulations confirming absolute concentration robustness in DOC systems}
\label{appendix:validation:ACR}

In this section, we validate our theoretical results for the absolute concentration robustness property in the DOC systems through simulations. As summarized in Table \ref{tab:summary} of the main manuscript, the positive and negative DOC systems do not exhibit ACR, consistent with the behavior shown in Figure \ref{fig:DOCGL}. For P-null systems, the table indicates that all species except $A_1$ exhibit ACR, while for Q-null systems, only species $A_1$ exhibits ACR. These findings are in agreement with the simulation results shown in Figures \ref{fig:DOCP0} and \ref{fig:DOCQ0}, respectively.

\begin{figure}[H]
    \begin{center}
    \includegraphics[width=16cm,height=14cm,keepaspectratio]{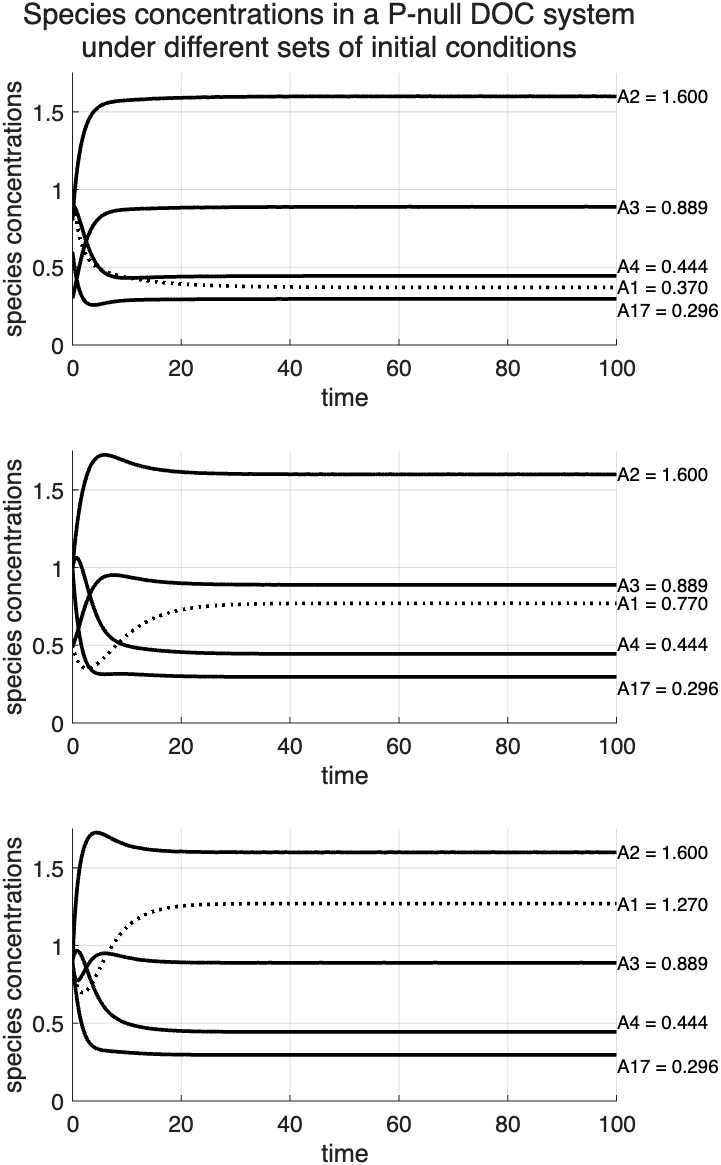}
    \captionsetup{font=scriptsize}
        \caption{Time evolution of species concentrations in a P-null DOC system simulated under three different sets of initial conditions. In this system, $p_1 = p_2$. Parameters for the rate constants: $k_1 = 0.5$, $k_2 = 0.8$, $k_3 = 0.5$, $k_4 = 0.7$, $k_5 = 0.4$, $k_6 = 0.6$, and $k_7 = 0.2$, and kinetic orders: $p_1 = 1.0$, $q_1 = 1.5$, $p_2 = 1.0$, and $q_2 = 0.5$ were used. The upper, middle, and lower subplots represent the system behavior under the following initial concentrations for $[A_1,\ A_2,\ A_3,\ A_4,\ A_{17}]$: 
    (upper) $[1.0,\ 0.8,\ 0.3,\ 0.9,\ 0.6]$, 
    (middle) $[0.5,\ 1.0,\ 0.5,\ 1.0,\ 1.0]$, and 
    (lower) $[0.9,\ 0.9,\ 0.9,\ 0.9,\ 0.9]$, respectively. The resulting steady-state concentrations verify that indeed the system exhibits absolute concentration robustness in all species except $A_1$.}
    \label{fig:DOCP0}
    \end{center}
    \end{figure}

\begin{figure}[H]
    \begin{center}
    \includegraphics[width=16cm,height=14cm,keepaspectratio]{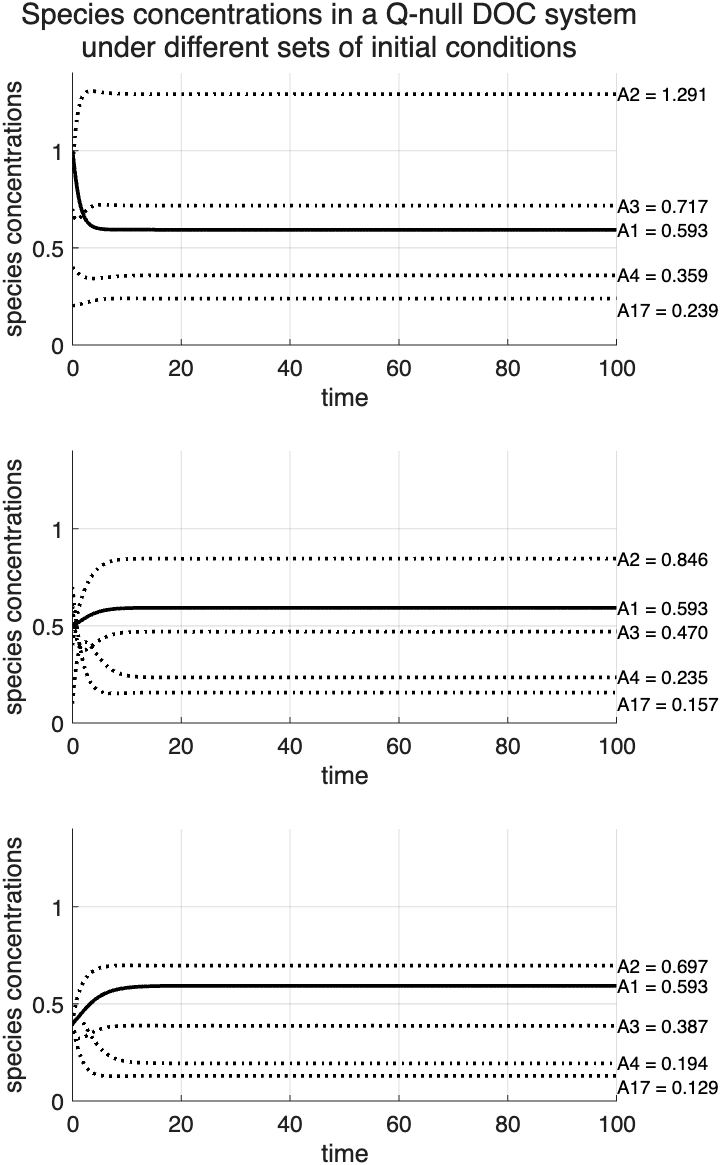}
    \captionsetup{font=scriptsize}
    \caption{Time evolution of species concentrations in a Q-null DOC system simulated under three different sets of initial conditions. In this system, $q_1 = q_2$. Parameters for the rate constants: $k_1 = 0.5$, $k_2 = 0.8$, $k_3 = 0.5$, $k_4 = 0.7$, $k_5 = 0.4$, $k_6 = 0.6$, and $k_7 = 0.2$, and kinetic orders: $p_1 = 0.5$, $q_1 = 1.5$, $p_2 = 1.4$, and $q_2 = 1.5$ were used. The upper, middle, and lower subplots represent the system behavior under the following initial concentrations for $[A_1,\ A_2,\ A_3,\ A_4,\ A_{17}]$: 
    (upper) $[1.0,\ 0.9,\ 0.7,\ 0.4,\ 0.2]$, 
    (middle) $[0.5,\ 0.4,\ 0.6,\ 0.1,\ 0.7]$, and 
    (lower) $[0.4,\ 0.4,\ 0.4,\ 0.4,\ 0.4]$, respectively. The observed steady-state concentrations reveal that the system exhibits absolute concentration robustness only in species $A_1$.}
    \label{fig:DOCQ0}
    \end{center}
    \end{figure}

\begin{figure}[H]
    \begin{center}
    \includegraphics[width=16cm,height=14cm,keepaspectratio]{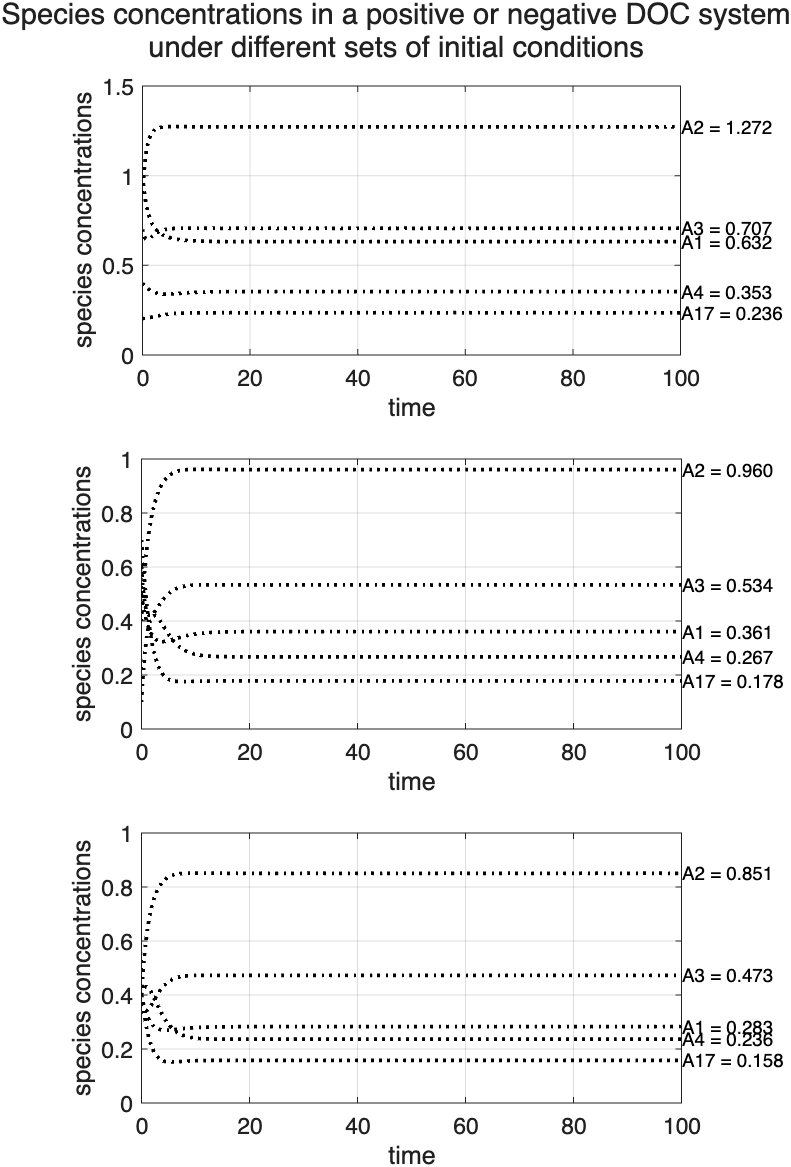}
    \captionsetup{font=scriptsize}
    \caption{Time evolution of species concentrations in a positive or negative DOC system simulated under three different sets of initial conditions. In this system, $p_1 \ne p_2$ and $q_1 \ne q_2$. Parameters for the rate constants: $k_1 = 0.5$, $k_2 = 0.8$, $k_3 = 0.5$, $k_4 = 0.7$, $k_5 = 0.4$, $k_6 = 0.6$, and $k_7 = 0.2$, and kinetic orders: $p_1 = 0.7$, $q_1 = 1.5$, $p_2 = 1.2$, and $q_2 = 0.5$ were used. The upper, middle, and lower subplots represent the system behavior under the following initial concentrations for $[A_1,\ A_2,\ A_3,\ A_4,\ A_{17}]$: 
    (upper) $[1.0,\ 0.9,\ 0.7,\ 0.4,\ 0.2]$, 
    (middle) $[0.5,\ 0.4,\ 0.6,\ 0.1,\ 0.7]$, and 
    (lower) $[0.4,\ 0.4,\ 0.4,\ 0.4,\ 0.4]$, respectively. The distinct steady-state concentrations observed in each subplot indicate that the system does not exhibit absolute concentration robustness.}\label{fig:DOCGL}
    \end{center}
    \end{figure}

\section{Positive steady state parametri-zation of the DAC system}

Following the steps provided in Appendix \ref{details:parametrization}, we also obtain a positive steady state parametrization of the DAC system, which was not provided in \cite{fortun2024determining}.
For the case when $p_1 \ne p_2$, we have
\begin{align*}
    \begin{cases}
    a_1 = \left(\dfrac{k_1}{k_2}\right)^{\frac{1}{p_2-p_1}} \tau^{q_2-q_1} \\
    a_2 = \tau^{p_1-p_2} \\
    a_3 = \dfrac{k_3}{k_4} \tau^{p_1-p_2} \\
    a_4 = \dfrac{k_6}{k_5} \tau^{p_1-p_2} \\
    a_5 = \dfrac{k_6}{k_7} \tau^{p_1-p_2} \\
    \text{Free parameter:} \hspace{0.2cm} \tau > 0.
    \end{cases}
\end{align*}

On the other hand, for the case when $q_1 \neq q_2$, we have
\begin{align*}
    \begin{cases}
    a_1 = \tau^{q_2-q_1} \\
    a_2 = \left(\dfrac{k_1}{k_2}\right)^{\frac{1}{q_2-q_1}} \tau^{p_1-p_2} \\
    a_3 = \dfrac{k_3}{k_4} \left(\dfrac{k_1}{k_2}\right)^{\frac{1}{q_2-q_1}} \tau^{p_1-p_2} \\
    a_4 = \dfrac{k_6}{k_5} \left(\dfrac{k_1}{k_2}\right)^{\frac{1}{q_2-q_1}} \tau^{p_1-p_2} \\
    a_5 = \dfrac{k_6}{k_7} \left(\dfrac{k_1}{k_2}\right)^{\frac{1}{q_2-q_1}} \tau^{p_1-p_2} \\
    \text{Free parameter:} \hspace{0.2cm} \tau > 0.
    \end{cases}
\end{align*}

\section{Positive steady state parametri-zation of the integrated system}
\label{details:param-combi}
Again, we follow the same steps in Appendix \ref{details:parametrization} to obtain a positive steady state parametrization of the system where both DOC and DAC technologies are present. If $p_1\neq p_2,$ then we have $a_2 = \tau^{p_1 - p_2}$ where $\tau > 0$ is a free parameter. Thus, we have
\begin{align*}
    \begin{dcases}
    a_1 = \left(\dfrac{k_1}{k_2}\right)^{\frac{1}{p_2-p_1}} \tau^{q_2-q_1} \\
    a_2 = \tau^{p_1-p_2} \\
    a_3 = \dfrac{k_3}{k_4 + k_7} \tau^{p_1-p_2} \\
    a_4 = \dfrac{k_3k_7 + k_4k_8 + k_7k_8}{k_5(k_4 + k_7)} \tau^{p_1-p_2} \\
    a_5 = \dfrac{k_8}{k_9} \tau^{p_1-p_2} \\
    a_{17} = \dfrac{k_3k_7}{k_6(k_4 + k_7)}\tau^{p_1 - p_2}\\
    \text{Free parameter:} \hspace{0.2cm} \tau > 0.
    \end{dcases}
\end{align*}
Next, if $q_1\neq q_2,$ then $a_2 = \left(\dfrac{k_1}{k_2}\right)^{\frac{1}{q_2-q_1}} \tau^{p_1-p_2}$, and so

\begin{align*}
    \begin{dcases}
    a_1 = \tau^{q_2 - q_1} \\
    a_2 = \left(\dfrac{k_1}{k_2}\right)^{\frac{1}{q_2-q_1}} \tau^{p_1-p_2} \\
    a_3 = \dfrac{k_3}{(k_4+k_7)}\left(\dfrac{k_1}{k_2}\right)^{\frac{1}{q_2-q_1}} \tau^{p_1-p_2}\\
    a_4 = \dfrac{k_3k_7 + k_4k_8 + k_7k_8}{k_5(k_4+k_7)}\left(\dfrac{k_1}{k_2}\right)^{\frac{1}{q_2-q_1}} \tau^{p_1-p_2} \\
    a_5 = \dfrac{k_8}{k_9}\left(\dfrac{k_1}{k_2}\right)^{\frac{1}{q_2-q_1}} \tau^{p_1-p_2}\\
    a_{17} = \dfrac{k_3k_7}{k_6(k_4+k_7)}\left(\dfrac{k_1}{k_2}\right)^{\frac{1}{q_2-q_1}} \tau^{p_1-p_2}\\
    \text{Free parameter:} \hspace{0.2cm} \tau > 0.
    \end{dcases}
\end{align*}

\section{Conditions of multistationarity for the integrated system}
\label{details:multi-combi}

We investigate the multistationary of the integrated system using the same arguments in Appendix \ref{details:multistationarity}.

Indeed, for positive systems, i.e. whenever $R > 0,$ we have multistationarity.

To determine some sufficient conditions for the direct ocean capture model to admit multiple steady states, we utilize a result by M\"uller and Regensburger \cite{regensburger}.

The Theorem \ref{muller} tells us that for weakly reversible generalized mass action systems, a sufficient condition for the system to be multistationary is the existence of a non-trivial vector whose sign pattern is the same as that of the stoichiometric subspace $S$ and the orthogonal complement of kinetic flux subspace $\tilde S.$

First, we solve for the sign pattern of $\tilde S.$ Note that $\tilde{S} = \text{Im}\,(\tilde{Y} \cdot I_a)$ where
\begin{center}
    $\tilde{Y}=\begin{blockarray}{ccccccccc}
        & \matindex{$A_1+2A_2$} & \matindex{$2A_1+A_2$} & \matindex{$A_2$} & \matindex{$A_3$} & \matindex{$A_4$} & \matindex{$A_5$} & \matindex{$A_{17}$} \\
        \begin{block}{c[cccccccc]}
        \matindex{$A_1$} & p_1 & p_2 & 0 & 0 & 0 & 0 & 0 \\
        \matindex{$A_2$} & q_1 & q_2 & 1 & 0 & 0 & 0 & 0 \\
        \matindex{$A_3$} & 0 & 0 & 0 & 1 & 0 & 0 & 0 \\
        \matindex{$A_4$} & 0 & 0 & 0 & 0 & 1 & 0 & 0 \\
        \matindex{$A_5$} & 0 & 0 & 0 & 0 & 0 & 1 & 0 \\
        \matindex{$A_{17}$} & 0 & 0 & 0 & 0 & 0 & 0 & 1 \\
        \end{block}
        \end{blockarray}$
\end{center}
and
\begin{center}
    $I_a=\begin{blockarray}{cccccccccc}
        & \matindex{$R_1$} & \matindex{$R_2$} & \matindex{$R_3$} & \matindex{$R_4$} & \matindex{$R_5$} & \matindex{$R_6$} & \matindex{$R_7$} & \matindex{$R_8$} & \matindex{$R_9$} \\
        \begin{block}{c[ccccccccc]}
        \matindex{$A_1+2A_2$} & -1 & 1 & 0 & 0 & 0 & 0 & 0 & 0 & 0 \\
        \matindex{$2A_1+A_2$} & 1 & -1 & 0 & 0 & 0 & 0 & 0 & 0 & 0 \\
        \matindex{$A_2$} & 0 & 0 & -1 & 1 & 1 & 0 & 0 & -1 & 0 \\
        \matindex{$A_3$} & 0 & 0 & 1 & -1 & 0 & 0 & -1 & 0 & 0 \\
        \matindex{$A_4$} & 0 & 0 & 0 & 0 & -1 & 1 & 0 & 0 & 1 \\
        \matindex{$A_5$} & 0 & 0 & 0 & 0 & 0 & 0 & 0 & 1 & -1 \\
        \matindex{$A_{17}$} & 0 & 0 & 0 & 0 & 0 & -1 & 1 & 0 & 0 \\
        \end{block}
        \end{blockarray}$.
\end{center} Here, the $\tilde Y$ matrix is defined using the kinetic order vectors of the system (see \cite{regensburger}) and $I_a$ is the incidence matrix of the network. Hence, 
$$\tilde{Y} \cdot I_a = \begin{bmatrix} 
p_2-p_1 & p_1-p_2 & 0 & 0 & 0 & 0 & 0 & 0 & 0 \\
q_2-q_1 & q_1-q_2 & -1 & 1 & 1 & 0 & 0 & -1 & 0 \\
0 & 0 & 1 & -1 & 0 & 0 & -1 & 0 & 0 \\
0 & 0 & 0 & 0 & -1 & 1 & 0 & 0 & 1 \\
0 & 0 & 0 & 0 & 0 & 0 & 0 & 1 & -1 \\
0 & 0 & 0 & 0 & 0 & -1 & 1 & 0 & 0
\end{bmatrix}$$
$$\Rightarrow \tilde{S} = \text{Im} \hspace{0.1cm} (\tilde{Y} \cdot I_a) = \sf{span} \hspace{0.1cm} \left \{\begin{bmatrix} p_2-p_1 \\ q_2-q_1 \\ 0 \\ 0 \\ 0 \\ 0 \end{bmatrix}, \begin{bmatrix} 0 \\ -1 \\ 1 \\ 0 \\ 0 \\ 0 \end{bmatrix}, \begin{bmatrix} 0 \\ 1 \\ 0 \\ -1 \\ 0 \\ 0 \end{bmatrix}, \begin{bmatrix} 0 \\ 0 \\ 0 \\ 1 \\ 0 \\ -1 \end{bmatrix}, \begin{bmatrix} 0 \\ -1 \\ 0 \\ 0 \\ 1 \\ 0 \end{bmatrix} \right \}.$$
The orthogonal complement $(\tilde{S})^{\perp}$ of $\tilde{S}$ is given by
$$(\tilde{S})^{\perp} = \rm{span} \left \{\begin{bmatrix} \frac{q_1-q_2}{p_2-p_1} \\ 1 \\ 1 \\ 1 \\ 1 \\ 1 \end{bmatrix} \right \} = \rm{span} 
 \left \{\begin{bmatrix} -Q \\ 1 \\ 1 \\ 1 \\ 1 \\ 1 \end{bmatrix} \right \} = \rm{span}  \left \{\begin{bmatrix} -1 \\ {R} \\ {R} \\ {R} \\ {R} \\ {R} \end{bmatrix} \right \}$$
where ${R} = \dfrac{p_2-p_1}{q_2-q_1}$ and ${Q} = \dfrac{q_2-q_1}{p_2-p_1},$ as defined.

We now investigate the multiplicity of steady states for positive ($R > 0$), negative ($R < 0)$, $P$-null ($R = 0$ and defined), and $Q$-null ($Q = 0$ and defined) systems. 

First, for positive integrated systems, i.e., $R > 0\, (Q > 0),$ we have \[{\rm sign}(\tilde S^\perp) =  \left \{\begin{bmatrix} - \\ + \\ + \\ + \\ + \\ + \end{bmatrix}, \begin{bmatrix} + \\ - \\ - \\ - \\ - \\ - \end{bmatrix} \right \}.\] 

Indeed, if we let $x$ be in the stoichiometric subspace $S$ given by \[S = \text{span }\left \{\begin{bmatrix} 1 \\ -1 \\ 0 \\ 0 \\ 0 \\ 0 \end{bmatrix}, \begin{bmatrix} 0 \\ -1 \\ 1 \\ 0 \\ 0 \\ 0 \end{bmatrix}, \begin{bmatrix} 0 \\ 1 \\ 0 \\ -1 \\ 0 \\ 0 \end{bmatrix}, \begin{bmatrix} 0 \\ 0 \\ 0 \\ 1 \\ 0 \\ -1 \end{bmatrix}, \begin{bmatrix} 0 \\ -1 \\ 0 \\ 0 \\ 1 \\ 0 \end{bmatrix} \right \},\] then 
\begin{center}
\resizebox{\textwidth}{!}{%
$
x = a_1 \begin{bmatrix} 1 \\ -1 \\ 0 \\ 0 \\ 0 \\ 0 \end{bmatrix}
  + a_2 \begin{bmatrix} 0 \\ -1 \\ 1 \\ 0 \\ 0 \\ 0 \end{bmatrix}
  + a_3 \begin{bmatrix} 0 \\ 1 \\ 0 \\ -1 \\ 0 \\ 0 \end{bmatrix}
  + a_4 \begin{bmatrix} 0 \\ 0 \\ 0 \\ 1 \\ 0 \\ -1 \end{bmatrix}
  + a_5 \begin{bmatrix} 0 \\ -1 \\ 0 \\ 0 \\ 1 \\ 0 \end{bmatrix}
  = \begin{bmatrix}
    a_1 \\
    -a_1 - a_2 + a_3 - a_5 \\
    a_2 \\
    -a_3 + a_4 \\
    a_5 \\
    -a_4
\end{bmatrix}.
$
}
\end{center}
We can then choose $a_1 > 0, a_2 < 0, a_3 > a_4 > 0, a_5 < 0$ such that $0 < -a_2 + a_3 -a_5 < a_1$ so that we have \[{\rm sign}(x) = \begin{bmatrix} + \\ - \\ - \\ - \\ - \\ - \end{bmatrix} \in {\rm sign}(\tilde S^\perp)\] and thus ${\rm sign}(x) \cap {\rm sign}(\tilde{S})^{\perp} \neq \{0\}.$ Therefore, by Theorem \ref{muller}, any positive integrated system is multistationary.

Now, for negative integrated systems, we cannot utilize Theorem \ref{muller} to conclude monostationarity. Because of this, we employ a different criterion to conclude when the system is monostationary. The following computational method introduced by Wiuf and Feliu \cite{wiuf, interacting} reveals network injectivity for a specific subset of the collection of negative integrated systems.

If we have a negative integrated system, i.e. we investigate network injectivity to assess multiplicity of steady states using Theorem \ref{powerlawinjectivity}.
For the integrated system, we similarly solve for $M^*$ as previously described and find that 

\[
\begin{aligned}
\det M^* = &\phantom{+}\,\,\, \underline{p_1}k_1 k_2 k_3 k_4 k_5 z_1 z_3 z_5 z_7 z_9 
+ \underline{p_1} k_1 k_2 k_3 k_4 k_6 z_1 z_4 z_5 z_6 z_8 \\
& -  \underline{p_2} k_1 k_2 k_3 k_4 k_5 z_2 z_3 z_5 z_7 z_9 
- \underline{p_2} k_1 k_2 k_3 k_4 k_6 z_2 z_4 z_5 z_6 z_8 \\
& + \underline{p_1} k_1 k_2 k_4 k_5 k_6  z_1 z_3 z_5 z_6 z_9 
+ \underline{p_1} k_1 k_2 k_3 k_4 k_6  z_1 z_5 z_6 z_7 z_8 \\
& + \underline{p_1} k_1 k_2 k_3 k_5 k_6 z_1 z_3 z_6 z_7 z_9 
- \underline{p_2} k_1 k_2 k_4 k_5 k_6 z_2 z_3 z_5 z_6 z_9 \\
& + \underline{p_1} k_1 k_3 k_4 k_5 k_6 z_1 z_4 z_5 z_6 z_9 
- \underline{p_2} k_1 k_2 k_3 k_4 k_6 z_2 z_5 z_6 z_7 z_8 \\
& + \underline{p_1} k_1 k_2 k_3 k_5 k_6 z_1 z_4 z_6 z_8 z_9 
- \underline{p_2} k_1 k_2 k_3 k_5 k_6 z_2 z_3 z_6 z_7 z_9 \\
& - \underline{p_2} k_1 k_3 k_4 k_5 k_6 z_2 z_4 z_5 z_6 z_9 
- \underline{p_2} k_1 k_2 k_3 k_5 k_6 z_2 z_4 z_6 z_8 z_9 \\
& + \underline{p_1} k_1 k_3 k_4 k_5 k_6 z_1 z_5 z_6 z_7 z_9 
+ \underline{p_1} k_1 k_2 k_3 k_5 k_6 z_1 z_6 z_7 z_8 z_9 \\
& - \underline{p_2} k_1 k_3 k_4 k_5 k_6 z_2 z_5 z_6 z_7 z_9 
- \underline{p_2} k_1 k_2 k_3 k_5 k_6 z_2 z_6 z_7 z_8 z_9 \\
& - \underline{q_1} k_2 k_3 k_4 k_5 k_6 z_1 z_4 z_5 z_6 z_9 
+ \underline{q_2} k_2 k_3 k_4 k_5 k_6 z_2 z_4 z_5 z_6 z_9 \\
& - \underline{q_1} k_2 k_3 k_4 k_5 k_6 z_1 z_5 z_6 z_7 z_9 
+ \underline{q_2} k_2 k_3 k_4 k_5 k_6 z_2 z_5 z_6 z_7 z_9.
\end{aligned}
\]
Similar to the steps in Appendix \ref{details:multistationarity}, all the terms of the determinant of $M^*$ are positive whenever $p_1 > 0, p_2 < 0, q_1 < 0,$ and $q_2 > 0.$ Similarly, all the terms are negative whenever $p_1 < 0, p_2 > 0, q_1 > 0,$ and $q_2 < 0.$ By Theorem \ref{powerlawinjectivity}, the integrated system is monostationary if either (i) $p_1, q_2 > 0$ and $p_2, q_1 < 0$ or (ii) $p_1, q_2 < 0$ and $p_2, q_1 > 0$ holds.

Finally, for null systems, we investigate the induced ODEs of the system and find that the integrated system follows the conservation law \[A'_1(t) + A'_2(t) + A'_3(t) + A'_4(t) + A'_5(t) + A'_{17}(t) = 0.\] Following the same arguments in Theorems \ref{thm-monoP0} and \ref{thm-monoQ0}, we can conclude that the integrated system whenever $p_1 = p_2$ or $q_1 = q_2,$ but not both, i.e. $P$-null or $Q$-null, admits a unique positive steady state.

 We notice that this is the same for the case of DOC-only systems.

\section{Carbon reduction for the DAC system}

\begin{proposition}
    \label{thm:DAC-reduction-suff}
    Let $A_i^\circ, A_i^*, S^\circ$ be as defined in the previous proposition of a DAC system. Let $m'$ be the minimum of ${\rm pr}_2$ and $M''$ be the maximum of ${\rm pr}_1 + {\rm pr}_3 + {\rm pr}_4 + {\rm pr}_5$ on $S^\circ.$ Then $A_2^* < A_2^\circ$ whenever \[1 + \frac{M''}{m'} < \left(\frac{k_1}{k_2}\right)^{\frac{1}{p_2 - p_1}}(m')^{-Q} + \frac{k_3k_5k_7 + k_6k_4(k_5 + k_7)}{k_4k_5k_7}\quad\text{if }p_1\neq p_2\] or \[1 + \frac{M''}{m'} < \left(\frac{k_1}{k_2}\right)^{\frac{1}{q_2 - q_1}}(m')^{-R} + \frac{k_3k_5k_7 + k_6k_4(k_5 + k_7)}{k_4k_5k_7}\quad\text{if }q_1\neq q_2.\]
    \begin{proof}
        Note that $A_2^* < A_2^\circ$ if and only if $1 < A_2^*/A_2^\circ.$ Equivalently, we have \[1 + \frac{({\rm pr}_1 + {\rm pr}_3 + {\rm pr}_4 + {\rm pr}_5)(A)}{A_2^*} < \frac{A_2^\circ}{A_2^*} + \frac{({\rm pr}_1 + {\rm pr}_3 + {\rm pr}_4 + {\rm pr}_5)(A)}{A_2^*}.\] Since $({\rm pr}_1 + {\rm pr}_3 + {\rm pr}_4 + {\rm pr}_5)(A) \leq M''$ and $1/A_2^* \leq 1/m'',$ we have \[A_2^* < A_2^\circ \iff 1 + \frac{M''}{m'} < \frac{A_2^\circ + ({\rm pr}_1 + {\rm pr}_3 + {\rm pr}_4 + {\rm pr}_5)(A)}{A_2^*}.\] Denote the right hand side of the above inequality as $\mathrm{SUM}^*_{(2)}.$ We establish a lower bound for $\mathrm{SUM}^*_{(2)}$ using the steady state parametrization of the system whenever $p_1\neq p_2$ with \begin{align*}&\left(\frac{k_1}{k_2}\right)^{1/p_2 - p_1}(m')^\frac{q_2 - q_1}{p_1 - p_2} + \frac{k_3}{k_4} + \frac{k_6}{k_5} + \frac{k_6}{k_7} \\
        &= \left(\frac{k_1}{k_2}\right)^{\frac{1}{p_2 - p_1}}(m')^{-Q} + \frac{k_3k_5k_7 + k_6k_4(k_5 + k_7)}{k_4k_5k_7}\leq \mathrm{SUM}^*_{(2)}.\end{align*} Similarly, if $q_1\neq q_2,$ we have \begin{align*}&\left(\frac{k_1}{k_2}\right)^{\frac{1}{q_2 - q_1}}(m')^{\frac{p_1 - p_2}{q_2 - q_1}} + \frac{k_3}{k_4} + \frac{k_6}{k_5} + \frac{k_6}{k_7} \\
        &= \left(\frac{k_1}{k_2}\right)^{\frac{1}{q_2 - q_1}}(m')^{-R} + \frac{k_3k_5k_7 + k_6k_4(k_5 + k_7)}{k_4k_5k_7}\leq \mathrm{SUM}^*_{(2)}.\end{align*} Then, the right hand side of the equivalence for $A_2^* < A_2^\circ$ is satisfied whenever \[1 + \frac{M''}{m'} < \left(\frac{k_1}{k_2}\right)^{\frac{1}{p_2 - p_1}}(m')^{-Q} + \frac{k_3k_5k_7 + k_6k_4(k_5 + k_7)}{k_4k_5k_7}\quad\text{if }p_1\neq p_2\] or \[1 + \frac{M''}{m'} < \left(\frac{k_1}{k_2}\right)^{\frac{1}{q_2 - q_1}}(m')^{-R} + \frac{k_3k_5k_7 + k_6k_4(k_5 + k_7)}{k_4k_5k_7}\quad\text{if }q_1\neq q_2.\]
    \end{proof}
\end{proposition} Note that Proposition \ref{thm:DAC-reduction-suff} implies that for null DAC systems, i.e. either $Q = 0$ or $R = 0,$ a sufficient condition for carbon reduction in the atmosphere is given by \[1 + \frac{M''}{m'} <\frac{k_3k_5k_7 + k_6k_4(k_5 + k_7)}{k_4k_5k_7}. \]


\singlespacing

\end{document}